\documentclass[a4paper,english]{extarticle}
\usepackage[utf8]{luainputenc}
\usepackage{color}
\usepackage{babel}
\usepackage{amsthm}
\usepackage{amsmath}
\usepackage{amssymb}
\usepackage{esint}
\PassOptionsToPackage{normalem}{ulem}
\usepackage{ulem}
\usepackage[unicode=true,pdfusetitle,
 bookmarks=true,bookmarksnumbered=false,bookmarksopen=false,
 breaklinks=false,pdfborder={0 0 1},backref=false,colorlinks=false]
 {hyperref}

\makeatletter

\pdfpageheight\paperheight
\pdfpagewidth\paperwidth

\numberwithin{equation}{section}
\numberwithin{figure}{section}
\theoremstyle{plain}
\newtheorem{thm}{\protect\theoremname}[section]
\theoremstyle{remark}
\newtheorem{rem}[thm]{\protect\remarkname}
\ifx\proof\undefined
\newenvironment{proof}[1][\protect\proofname]{\par
\normalfont\topsep6\p@\@plus6\p@\relax
\trivlist
\itemindent\parindent
\item[\hskip\labelsep\scshape #1]\ignorespaces
}{%
\endtrivlist\@endpefalse
}
\providecommand{\proofname}{Proof}
\fi
\theoremstyle{plain}
\newtheorem{lem}[thm]{\protect\lemmaname}
\theoremstyle{plain}
\newtheorem{prop}[thm]{\protect\propositionname}
\theoremstyle{definition}
\newtheorem{defn}[thm]{\protect\definitionname}

\usepackage{cancel}
\numberwithin{equation}{section}
\date{}

\makeatother

\providecommand{\definitionname}{Definition}
\providecommand{\lemmaname}{Lemma}
\providecommand{\propositionname}{Proposition}
\providecommand{\remarkname}{Remark}
\providecommand{\theoremname}{Theorem}

\begin{document}

\global\long\def\apa{\alpha}
\global\long\def\ba{\beta}
\global\long\def\ga{\gamma}
\global\long\def\Ga{\Gamma}
\global\long\def\da{\delta}
\global\long\def\Da{\Delta}
\global\long\def\ta{\theta}
\global\long\def\Ta{\Theta}
\global\long\def\la{\lambda}

\global\long\def\La{\Lambda}
\global\long\def\oa{\omega}
\global\long\def\Oa{\Omega}
\global\long\def\vn{\varepsilon}
\global\long\def\vi{\varphi}
\global\long\def\en{\epsilon}


\global\long\def\pl{\partial}
\global\long\def\na{\nabla}
\global\long\def\sto{\rightsquigarrow}
\global\long\def\wto{\rightharpoonup}
\global\long\def\To{\Rightarrow}
\global\long\def\Lra{\Leftrightarrow}

\global\long\def\fc#1#2{\frac{#1}{#2}}
\global\long\def\st#1{\sqrt{#1}}
\global\long\def\ol#1{\overline{#1}}
\global\long\def\ul#1{\underline{#1}}
\global\long\def\wt#1{\widetilde{#1}}

\global\long\def\wh#1{\widehat{#1}}


\global\long\def\mf#1{\mathbf{#1}}
\global\long\def\mb#1{\mathbb{#1}}
\global\long\def\ml#1{\mathcal{#1}}
\global\long\def\mk#1{\mathfrak{#1}}
\global\long\def\mm#1{\mathrm{#1}}
\global\long\def\bs#1{\boldsymbol{#1}}

\title{The enclosure method for the anisotropic Maxwell system}

\author{Rulin Kuan; Yi-Hsuan Lin; Mourad Sini}
\maketitle
\begin{abstract}
We develop an enclosure-type reconstruction scheme to identify penetrable
and impenetrable obstacles in electromagnetic field with anisotropic
medium in $\mathbb{R}^{3}$. The main difficulty in treating this
problem lies in the fact that there are so far no complex geometrical
optics solutions available for the Maxwell's equation with anisotropic
medium in $\mathbb{R}^{3}$. Instead, we derive and use another type
of special solutions called oscillating-decaying solutions. To justify
this scheme, we use Meyers' $L^{p}$ estimate, for the Maxwell system,
to compare the integrals coming from oscillating-decaying solutions
and those from the reflected solutions. 
\end{abstract}

\section{Introduction and statement of the results}

\textcolor{black}{Let $\Oa$ be a bounded $C^{\infty}$-smooth domain
in $\mb R^{3}$ with connected complement $\mb R^{3}\setminus\ol{\Oa}$
and $D$ be a subset of $\Oa$ with Lipschitz boundary. We are concerned
with the electromagnetic wave propagation in an anisotropic medium
in $\mathbb{R}^{3}$ with the electric permittivity $\epsilon=(\epsilon_{ij}(x))$
a $3\times3$ positive definite matrix and $\epsilon(x)=\epsilon_{0}(x)$
in $\Omega\backslash\bar{D}$. We also assume that $\epsilon(x)=\epsilon_{0}(x)-\epsilon_{D}(x)\chi_{D}(x)$
with $\epsilon_{0}\in C^{\infty}(\Omega)$ a positive definite $3\times3$
symmetric matrix and $\epsilon_{D}(x)$ is a positive $3\times3$
symmetric matrix and }$\mu$ a smooth scalar function defined on \textcolor{black}{$\Omega$
such that there }exist $\mu_{c}>0$ and $\epsilon_{c}>0$ verifying
\begin{equation}
\mu(x)\geq\mu_{c}>0\mbox{ and }\sum_{i.j=1}^{3}\epsilon_{ij}(x)\xi_{i}\xi_{j}\geq\epsilon_{c}|\xi|^{2}\mbox{ }\forall\xi\in\mathbb{R}^{3},\mbox{ }\forall x\in\Omega.\label{eq:0.0}
\end{equation}
\textcolor{black}{If we denote by $E$ and $H$ the electric and the
magnetic fields respectively, then the electromagnetic wave propagation
by a penetrable obstacle problem reads as 
\begin{equation}
\begin{cases}
\nabla\times E-ik\mu H=0 & \mbox{ in }\Omega,\\
\nabla\times H+ik\epsilon E=0 & \mbox{ in }\Omega,\\
\nu\times E=f & \mbox{ on }\partial\Omega,
\end{cases}\label{eq:0.1}
\end{equation}
with $\epsilon=\epsilon_{0}-\epsilon_{D}\chi_{D}$, and the one by
the impenetrable obstacle as 
\begin{equation}
\begin{cases}
\nabla\times E-ik\mu H=0 & \mbox{ in }\Omega\backslash\bar{D},\\
\nabla\times H+ik\epsilon E=0 & \mbox{ in }\Omega\backslash\bar{D},\\
\nu\times E=f & \mbox{ on }\partial\Omega,\\
\nu\times H=0 & \mbox{ on }\partial D,
\end{cases}\label{eq:0.2}
\end{equation}
where $\nu$ is the unit outer normal vector on $\partial\Omega\cup\partial D$
and $k>0$ is the wave number. In this paper, we assume that $k$
is not an eigenvalue for (\ref{eq:0.1}) and (\ref{eq:0.2}).}\\
 \textbf{Impedance Map}: We define the impedance map $\Lambda_{D}:TH^{-\frac{1}{2}}(\partial\Omega)\to TH^{-\frac{1}{2}}(\partial\Omega)$
by 
\[
\Lambda_{D}(\nu\times H|_{\partial\Omega})=(\nu\times E|_{\partial\Omega}),
\]
where $TH^{-\frac{1}{2}}(\partial\Omega):=\{f\in H^{-\frac{1}{2}}(\partial\Omega)|\nu\cdot f=0\}$
and $\times$ is the standard cross product in $\mathbb{R}^{3}$.
We denote by $\Lambda_{\emptyset}$ the impedance map for the domain
without an obstacle.

Consider the anisotropic Maxwell system 
\begin{equation}
\begin{cases}
\nabla\times E-ik\mu H=0 & \mbox{ in }\Omega,\\
\nabla\times H+ik\epsilon E=0 & \mbox{ in }\Omega,
\end{cases}\label{eq:main}
\end{equation}
where $\mu$ and $\epsilon$ satisfy (\ref{eq:0.0}). We are interested
in the question reconstructing the shape of $D$ using the impedance
map $\Lambda_{D}$. This geometrical inverse problem is quite well
studied in the literature see \cite{isakov2006inverse} and several
methods have been proposed to solve it. In this paper, we focus on
one of these method, called the enclosure method, which is initiated
by Ikehata, see for examples \cite{I1999,I1998}, and developed by
many researchers \cite{KS2014,K2012,NY2007,SW2006,SY2012reconstruction,UW2008},
\cite{kar2014farfield,SY2012reconstruction} for the acoustic model,
\cite{kar2013elasticfarfield,K2012} for the Lam$\acute{\mathrm{e}}$
model and \cite{KS2014,Z2010} for the Maxwell model. The testing
functions used in \cite{KS2014,Z2010} are complex geometric optics
(CGO) solutions of the isotropic Maxwell's equation. The construction
of CGO solutions for isotropic inhomogeneous Maxwell's equations is
first proposed in \cite{OS1996}. After that, the authors in \cite{KSU2011}
also constructed CGO solutions for some special anisotropic Maxwell's
equations. However, there are not yet of CGO solutions for general
anisotropic Maxwell system. Besides, CGO solutions, another kind of
special solutions for anisotropic elliptic system was proposed for
substitution in \cite{NUW2005(ODS)} and \cite{NUW2006}. They are
called oscillating-decaying (OD) solutions. Inspired by \cite{OS1996}
and \cite{NUW2005(ODS)}, our idea is to reduce (\ref{eq:main}) to
an elliptic systems and then use the results in \cite{NUW2005(ODS)}
to construct oscillating-decaying type solutions to the anisotropic
Maxwell system. Precisely, we can decompose the equation (\ref{eq:main})
into two decoupled strongly elliptic systems. The main difference
between the construction of the oscillating-decaying solutions in
\cite{NUW2005(ODS)} and ours is about the higher derivatives of oscillating-decaying
solutions.

One of the main differences between the CGOs and the oscillating-decaying
solutions is that, roughly speaking, given a hyperplane, an oscillating
decaying solution is oscillating very rapidly along this plane and
decaying exponentially in the direction transversely to the same plane.
Oscillating-decaying solutions are special solutions with the phase
function having nonnegative imaginary part. In addition, these oscillating
decaying solutions are only defined on a half plane. To use them as
inputs for our detection algorithm, we need to extend them to the
whole domain $\Omega$. One way to do the extension is to use the
Runge approximation property for the anisotropic Maxwell's equation.
The Runge approximation property will help us to find a sequence of
approximated solutions which are defined on $\Omega$, satisfy (\ref{eq:main})
and their limit is the oscillating-decaying solution. Note that it
was first recognized by Lax \cite{lax1956stability} that the Runge
approximation property is a consequence of the weak unique continuation
property. In \cite{leis2013initial}, the authors already proved the
unique continuation property and based on it we derive the Runge approximation
property for the anisotropic Maxwell's equation.

To be more precise, let $\omega$ be a unit vector in $\mathbb{R}^{3}$,
denote $\Omega_{t}(\omega)=\Omega\cap\{x|x\cdot\omega>t\}$, $\Sigma_{t}(\omega)=\Omega\cap\{x|x\cdot\omega=t\}$
and set $(E_{t},H_{t})$ to be the oscillating-decaying solution for
the anisotropic Maxwell's equation in $\Omega_{t}(\omega)$.\\
 \textbf{Support function}: For $\rho\in\mathbb{S}^{2}$, we define
the support function of $D$ by $h_{D}(\rho)=\inf_{x\in D}x\cdot\rho$.

When $t=h_{D}(\rho)$, which means $\Sigma_{t}(\omega)$ touches $\partial D$,
we cannot apply the Runge approximation property to $(E_{t},H_{t})$
in $\Omega_{t}(\omega)$. Therefore, we need to enlarge the domain
$\Omega_{t}(\omega)$ such that the OD solutions exist and the Runge
approximation property works. Let $\eta$ be a positive real number,
denote $\Omega_{t-\eta}(\omega)$ and $\Sigma_{t-\eta}(\omega)$ and
note that $\Omega_{t-\eta}(\omega)\subset\Omega_{t}(\omega)$ $\forall\eta>0$.
We can find $(E_{t-\eta},H_{t-\eta})$ to be the OD solution in $\Omega_{t-\eta}(\omega)$.
By the Runge approximation property, there exists a sequence of functions
$\{(E_{\eta,\ell},H_{\eta,\ell})\}$ satisfying the Maxwell system
in $\Omega$ such that $(E_{\eta,\ell},H_{\eta,\ell})$ converges
to $(E_{t-\eta},H_{t-\eta})$ as $\ell\to\infty$ in $L^{2}(\Omega_{t-\eta}(\omega))$
and in $H(curl,D)$ by interior estimates since $D\Subset\Omega_{t-\eta}(\omega)$.
In addition we show that $(E_{t-\eta},H_{t-\eta})$ converges to $(E_{t},H_{t})$
in $H(curl,D)$ as $\eta\to0$.\textbf{ }Then we can define the indicator
function as follows.\\
 \textbf{Indicator function}: For $\rho\in S^{2}$, $\tau>0$ and
$t>0$ we define the indicator function 
\[
I_{\rho}(\tau,t):=\lim_{\eta\to0}\lim_{\ell\to\infty}I_{\rho}^{\eta,\ell}(\tau,t),
\]
where 
\[
I_{\rho}^{\eta,\ell}(\tau,t):=ik\tau\int_{\partial\Omega}(\nu\times H_{\eta,\ell})\cdot(\overline{(\Lambda_{D}-\Lambda_{\emptyset})(\nu\times H_{\eta,\ell})}\times\nu)dS.
\]
\textbf{Goal}: We want to characterize the convex hull of the obstacle
$D$ from the impedance map $\Lambda_{D}$.

The answer to this goal is the following theorem.
\begin{thm}
\label{main-thm} Let $\rho\in\mathbb{S}^{2}$. For the penetrable
(or impenetrable) obstacle case, we have the following characterization
of $h_{D}(\rho)$.
\[
\begin{cases}
\lim_{\tau\to\infty}|I_{\rho}(\tau,t)|=0\mbox{ when }t<h_{D}(\rho),\\
\liminf_{\tau\to\infty}|I_{\rho}(\tau,h_{D}(\rho))|>0,
\end{cases}
\]

\end{thm}
To prove Theorem \ref{main-thm}, for the penetrable obstacle case,
we need an appropriate $L^{p}$ estimate of the corresponding reflected
solution. We follow the idea in \cite{KS2014} to prove a global $L^{p}$
estimate for the curl of the solutions of the anisotropic Maxwell's
equation, for $p$ near 2 and $p\leq2$.

To prove Theorem \ref{main-thm}, in the impenetrable obstacle case,
we use layer potential arguments as in \cite{KS2014} coupled with
appropriate $L^{p}$ estimates. Precisely, first, we use the well-posedness
for an exterior isotropic Maxwell's system with the Silver-M$\ddot{\mathrm{u}}$ller
radiation condition and, in particular, the layer potential theory
to find a suitable estimate for the solution of this exterior problem.
Second, we decompose the reflected solution into two functions, one
satisfies the reflected Maxwell's equation with a zero boundary data,
the other satisfies the original anisotropic Maxwell's equation with
the same boundary conditions which come from the reflected equation.
For the first decomposed function, we use the $L^{p}$ estimates,
and for the second function, we will use the well-posedness, in $L^{2}$,
for the anisotropic Maxwell's system. Combining these two steps, we
derive the full estimate for the reflected solution in the impenetrable
obstacle case.

This paper is organized as follows. In the section 2, we give decompose
the anisotropic Maxwell system into two strongly elliptic systems.
In section 3, we use the elliptic systems derived in the section 2
to build the oscillating-decaying solutions for the Maxwell system.
Then, we give the Runge approximation for the anisotropic Maxwell
equation in section 4. In section 5, we prove the Theorem \ref{main-thm}
for both penetrable and impenetrable obstacle case. Finally, in the
last section, as an appendix, we provide some technical details which
we postponed in the main text and recall some useful estimates for
solutions of the Maxwell system. Before closing this introduction,
let us mention that in the whole text whenever we use the word smooth
it means $C^{\infty}$-smooth.

\section{Reduction to strongly elliptic systems }

Our goal is to construct the oscillating-decaying (OD) solution for
the following anisotropic time-harmonic Maxwell's system 
\begin{align}
\begin{cases}
\begin{aligned}\na\times E=ik\mu H\\
\na\times H=-ik\en E\\
\mbox{div}(\en E)=0\\
\mbox{div}(\mu H)=0,
\end{aligned}
\end{cases}.\label{ME}
\end{align}
where $E,H$ denote the electric and magnetic field intensity respectively,
and $\mu$ denotes the positive scalar permeability, $\en$ denotes
the permittivity, which is a real, symmetric, positive definite $3\times3$
matrix.

Inspired by \cite{OS1996}, the first step of constructing OD solutions
is to reduce (\ref{ME}) to a strongly elliptic system. In fact, we
reduce the anisotropic Maxwell's system (\ref{ME}) to two separate
strongly elliptic equations (\ref{eq:1.3}), while in \cite{OS1996}
the isotropic Maxwell's system is reduced to an elliptic (a single
Schrödinger) system with coupled zero-th order term. The following
theorem is our reduction result.
\begin{thm}
We set $E$ and $H$ of the following forms 
\begin{equation}
\begin{cases}
E=-\dfrac{i}{k}\en^{-1}\nabla\times(\mu^{-1}(\nabla\times B))-\en^{-1}(\nabla\times A)\\
H=\dfrac{i}{k}\mu^{-1}\nabla\times(\en^{-1}(\nabla\times A))-\mu^{-1}(\nabla\times B)
\end{cases}\label{eq:1.2}
\end{equation}
with $A,B$ satisfying the strongly elliptic systems 
\begin{equation}
\begin{cases}
\mu\nabla tr(M^{A}\nabla A)-\nabla\times(\en^{-1}(\nabla\times A))+k^{2}\mu A=0\\
\en\nabla tr(M^{B}\nabla B)-\nabla\times(\mu^{-1}(\nabla\times B))+k^{2}\en B=0
\end{cases},\label{eq:1.3}
\end{equation}
where $M^{A},M^{B}$ are introduced in Theorem \ref{Max-SE}, then
$E$ and $H$ satisfy (\ref{ME}).\end{thm}
\begin{rem}
Theorem 2.1 shows that, if we can find solutions of (\ref{eq:1.3}),
then we can find solutions of (\ref{ME}).\end{rem}
\begin{proof}
In this proof, we will show the process of the reduction. And the
proof that the systems (\ref{eq:1.3}) are strongly elliptic systems
will be postponed to Theorem \ref{Max-SE}.

As in \cite{OS1996}, we set the following two auxiliary functions
which are similar to what they used: 
\[
\Phi=\frac{i}{k}\mbox{div}(\en E)
\]
and 
\[
\Psi=\frac{i}{k}\mbox{div}(\mu H).
\]
Note that $\Phi$ and $\Psi$ are actually zero by the Maxwell's equation.
We consider the following first-order matrix differential operator
$P$ 
\[
P=\left(\begin{array}{cccc}
0 & \mbox{div}(\en(\cdot)) & 0 & 0\\
\mu^{-1}\na & 0 & \na\times & 0\\
0 & -\na\times & 0 & \en^{-1}\na\\
0 & 0 & \mbox{div}(\mu(\cdot)) & 0
\end{array}\right).
\]
Note that $P$ is a $8\times8$ matrix. Let 
\[
Y=\left(\begin{aligned}\Phi\\
E\\
H\\
\Psi
\end{aligned}
\right)
\]
Then the problem (\ref{ME}) can be rewritten as follows: 
\[
PY=-ikVY,
\]
where 
\[
V=\left(\begin{array}{cccc}
1 & 0 & 0 & 0\\
0 & \en & 0 & 0\\
0 & 0 & \mu & 0\\
0 & 0 & 0 & 1
\end{array}\right)
\]
Thus, the Maxwell's system (\ref{ME}) implies 
\begin{align}
(P+ikV)Y=0\mbox{ and }\Phi=\Psi=0.\label{PE}
\end{align}
It is easy to see that conversely (\ref{PE}) implies the Maxwell's
system, and hence they are equivalent.

The first idea of the reducing process is to construct a suitable
$\wt Q$, which can make $(P+ikV)\wt Q$ a ``good'' second-order
differential operator. Then, a solution $\mf X$ for the problem 
\begin{align}
(P+ikV)\wt QX=0\label{PQeq}
\end{align}
will give rise to a solution $Y=\wt QX$ for 
\[
(P+ikV)Y=0.
\]
Moreover, if we find the solution $X$ such that the first and the
last component of $Y=\wt QX$ are zero, then we obtain solutions for
the Maxwell's system.

We try the matrix differential operator $\wt Q=Q-ikI$, where 
\begin{equation}
Q=\left(\begin{array}{cccc}
0 & \mbox{div}(\en(\cdot)) & 0 & 0\\
\na & 0 & \en^{-1}(\na\times(\cdot)) & 0\\
0 & -\mu^{-1}(\na\times(\cdot)) & 0 & \na\\
0 & 0 & \mbox{div}(\mu(\cdot)) & 0
\end{array}\right).\label{eq:Q}
\end{equation}
Then 
\begin{align*}
 & (P+ikV)\wt Q\\
 & =(P+ikV)(Q-ikI)\\
 & =PQ-ikP+ikVQ+k^{2}V\\
 & =\left(\begin{array}{cccc}
\mbox{div}(\en\na) & 0 & 0 & 0\\
0 & L_{1} & 0 & 0\\
0 & 0 & L_{2} & 0\\
0 & 0 & 0 & \mbox{div}(\mu\na)
\end{array}\right)\\
 & +\left(\begin{array}{cccc}
0 & -ik\mbox{div}(\en(\cdot)) & 0 & 0\\
-ik\mu^{-1}\na & 0 & -ik\na\times & 0\\
0 & ik\na\times & 0 & -ik\en^{-1}\na\\
0 & 0 & -ik\mbox{div}(\mu(\cdot)) & 0
\end{array}\right)\\
 & +\left(\begin{array}{cccc}
0 & ik\mbox{div}(\en(\cdot)) & 0 & 0\\
ik\en\na & 0 & ik\na\times & 0\\
0 & -ik\na\times & 0 & ik\mu\na\\
0 & 0 & ik\mbox{div}(\mu(\cdot)) & 0
\end{array}\right)\\
 & +\left(\begin{array}{cccc}
k^{2} & 0 & 0 & 0\\
0 & k^{2}\en & 0 & 0\\
0 & 0 & k^{2}\mu & 0\\
0 & 0 & 0 & k^{2}
\end{array}\right)\\
 & =\left(\begin{array}{cccc}
\mbox{div}(\en\na)+k^{2} & 0 & 0 & 0\\
ik(\en-\mu^{-1})\na & L_{1}+k^{2}\en & 0 & 0\\
0 & 0 & L_{2}+k^{2}\mu & ik(\mu-\en^{-1})\na\\
0 & 0 & 0 & \mbox{div}(\mu\na)+k^{2}
\end{array}\right),
\end{align*}
where 
\begin{align}
 & L_{1}=\mu^{-1}\na(\mbox{div}(\en(\cdot)))-\na\times(\mu^{-1}(\na\times(\cdot)))\label{Li}\\
 & L_{2}=\en^{-1}\na(\mbox{div}(\mu(\cdot)))-\na\times(\en^{-1}(\na\times(\cdot))).
\end{align}

A prominent feature of the above operator is that it decomposes the
original eight-component system into two four-component systems. Precisely,
Set 
\[
X=\left(\begin{aligned}\varphi\\
e\\
h\\
\psi
\end{aligned}
\right),
\]
then (\ref{PQeq}) can be separated into two systems: 
\begin{align*}
\left\{ \begin{aligned}\mbox{div}(\en\na\vi)+k^{2}\vi=0\\
L_{1}e+k^{2}\en e+ik(\en-\mu^{-1})\na\vi=0.
\end{aligned}
\right.
\end{align*}
and 
\begin{align*}
\left\{ \begin{aligned}\mbox{div}(\mu\na\psi)+k^{2}\psi=0\\
L_{2}h+k^{2}\mu h+ik(\mu-\en^{-1})\na\psi=0.
\end{aligned}
\right.
\end{align*}
Moreover, 
\begin{align*}
Y & =\wt QX\\
 & =\left[\left(\begin{array}{cccc}
0 & \mbox{div}(\en(\cdot)) & 0 & 0\\
\na & 0 & \en^{-1}(\na\times(\cdot)) & 0\\
0 & -\mu^{-1}(\na\times(\cdot)) & 0 & \na\\
0 & 0 & \mbox{div}(\mu(\cdot)) & 0
\end{array}\right)-ikI\right]X\\
 & =\left(\begin{array}{l}
\mbox{div}(\en e)-ik\vi\\
\na\vi+\en^{-1}(\na\times h)-ike\\
-\mu^{-1}(\na\times\mf e)+\na\psi-ikh\\
\mbox{div}(\mu h)-ik\psi
\end{array}\right).
\end{align*}
Therefore, the problem of finding the solutions $X$ of 
\begin{align}
(P+ikV)\wt QX=0\mbox{ with the first and last component of }\wt QX\mbox{ being }0\label{PQ0}
\end{align}
is equivalent to the problem of finding solutions of the following
two separate systems: 
\begin{align}
\left\{ \begin{array}{l}
\mbox{div}(\en e)-ik\vi=0,\\
\mbox{div}(\en\na\vi)+k^{2}\vi=0,\\
\mu^{-1}\na(\mbox{div}(\en e))-\na\times(\mu^{-1}(\na\times e))+k^{2}\en e+ik(\en-\mu^{-1})\na\vi=0,
\end{array}\right.\label{un-e}
\end{align}
and 
\begin{align}
\left\{ \begin{array}{l}
\mbox{div}(\mu h)-ik\psi=0,\\
\mbox{div}(\mu\na\psi)+k^{2}\psi=0,\\
\en^{-1}\na(\mbox{div}(\mu h))-\na\times(\en^{-1}(\na\times h))+k^{2}\mu h+ik(\mu-\en^{-1})\na\psi=0.
\end{array}\right.\label{un-h}
\end{align}
Notice that if we set $e$ in the following form 
\begin{align}
e=-\fc ik(\na\vi+\en^{-1}(\na\times A)),\label{heform-e}
\end{align}
then the first equation of (\ref{un-e}) becomes the same as the second
one. For the third equation, we have 
\begin{align*}
 & \mu^{-1}\na\big(\mbox{div}(\en e)\big)-\na\times\big(\mu^{-1}(\na\times e)\big)+k^{2}\en e+ik(\en-\mu^{-1})\na\vi\\
 & =-\fc ik\mu^{-1}\na\left(\mbox{div}(\en\na\vi)\right)+\fc ik\na\times\bigg(\mu^{-1}\big[\na\times(\en^{-1}(\na\times A))\big]\bigg)\\
 & \qquad-ik\en\na\vi-ik(\na\times A)+ik\en\na\vi-\fc ik\mu^{-1}\na\Big(k^{2}\vi\Big)\\
 & =-\fc ik\mu^{-1}\na\left(\mbox{div}(\en\na\vi)+k^{2}\vi\right)+\fc ik\na\times\bigg(\mu^{-1}\big[\na\times(\en^{-1}(\na\times A))\big]\bigg)-ik\na\times A\\
 & =0+\fc ik\na\times\bigg(\mu^{-1}\big[\na\times(\en^{-1}(\na\times A))\big]\bigg)-ik\na\times A,
\end{align*}
by the second equation of (\ref{un-e}). Thus, by letting $e$ be
of the form (\ref{heform-e}), the system (\ref{un-e}) reduces to
\begin{align}
\left\{ \begin{array}{l}
\mbox{div}\big(\ga\na\vi\big)+k^{2}\vi=0,\\
\na\times\bigg(\mu^{-1}\big[\na\times(\en^{-1}(\na\times A))\big]-k^{2}A\bigg)=0.
\end{array}\right.\label{viA'}
\end{align}
Similarly, by letting 
\begin{align*}
h=-\fc ik(\na\psi+\mu^{-1}(\na\times B))
\end{align*}
for some vector field $B$, we can reduce (\ref{un-h}) to the following
system: 
\begin{align}
\left\{ \begin{array}{l}
\mbox{div}\big(\mu\na\psi\big)+k^{2}\psi=0,\\
\na\times\bigg(\en^{-1}\big[\na\times(\mu^{-1}(\na\times B))\big]-k^{2}B\bigg)=0.
\end{array}\right.\label{psB'}
\end{align}
To resume, if we can find solutions $\vi,A,\psi$ and $B$ of (\ref{viA'})
and (\ref{psB'}), we can find solutions of the problem (\ref{PQ0})
and therefore the original problem (\ref{ME}).

Now let us focus on (\ref{viA'}) and (\ref{psB'}). The goal is to
find special solutions (e.g. oscillating-decaying solutions) of (\ref{viA'})
and (\ref{psB'}). The idea of doing that is to subtract zero terms
of the form $\na\times\big(\na tr(M^{A}\na A)\big)$ and $\na\times\big(\na tr(M^{B}\na B)\big)$
from the second equations of (\ref{viA'}) and (\ref{psB'}) for some
matrices $M^{A},M^{B}$, so that they become $\na\times(\ml L^{A}A)=0$
and $\na\times(\ml L^{B}B)=0$ with $\ml L^{A}$ and $\ml L^{B}$
being strongly elliptic operators. Precisely, we want to find suitable
matrices $M^{A}$ and $M^{B}$ such that 
\begin{align}
\mu\na tr(M^{A}\na A)-\na\times\big(\en^{-1}(\na\times A)\big)+k^{2}\mu A=0\label{mA}
\end{align}
and 
\begin{align}
\en\na tr(M^{B}\na B)-\na\times\big(\mu^{-1}(\na\times B)\big)+k^{2}\en B=0\label{mB}
\end{align}
are strongly elliptic systems. In fact, by letting $M^{A}=m\mu^{-1}I$
and $M^{B}=m\mu^{-1}\en$, we can show that (\ref{mA}) and (\ref{mB})
are strong elliptic systems for arbitrary positive constant $m$.
The proof are given in Theorem \ref{Max-SE}. 
\end{proof}
To prove Theorem \ref{Max-SE}, we start with the following computational
lemma.
\begin{lem}
\label{id-curl} Let $M$ be a matrix-valued function with smooth
entries and $\mf F$ be a vector field. Then the $i$-th component
of the vector $\nabla\times\big(M(\nabla\times{\bf {F})\big)}$ is
given by 
\begin{align}
\left(\nabla\times\big(M(\nabla\times{\bf F})\big)\right){}_{i}=\sum_{j,k,\ell}\wt C_{ijk\ell}\pl_{j\ell}f_{k}+\wt R_{i},
\end{align}
where 
\begin{align*}
\wt C_{ijk\ell}=\da_{j\ell}M_{ki}+\da_{ik}M_{\ell j}-\da_{jk}M_{\ell i}-\da_{i\ell}M_{kj}+\big(\da_{i\ell}\da_{jk}-\da_{ik}\da_{j\ell}\big)tr(M),
\end{align*}
and $\wt R_{i}$ contains the lower order terms. Here, $\da_{ij}$
is the Kronecker delta, $M_{ij}$ is the $ij$-th entry of $M$, and
$\mf F=(f_{1},f_{2},f_{3})^{T}$. \end{lem}
\begin{proof}
We prove it by direct computations. For any vectors $\mf a,\mf b$,
letting $\mf c=\mf a\times\mf b$, we have 
\begin{align*}
c_{m}=\sum_{k\ell}\vn_{m\ell k}a_{\ell}b_{k},
\end{align*}
where $\mf a=(a_{1},a_{2},a_{3})^{T}$, $\mf b=(b_{1},b_{2},b_{3})^{T}$,
$\mf c=(c_{1},c_{2},c_{3})^{T}$ and $\vn_{m\ell k}$ denotes the
Levi-Civita symbol. Therefore, we obtain the $m$-th component of
$\na\times\mf F$: 
\begin{align*}
\bigg(\na\times\mf F\bigg)_{m}=\sum_{k\ell}\vn_{m\ell k}\pl_{\ell}f_{k}.
\end{align*}
Then, the $n$-th component of $M(\na\times\mf F)$ is 
\begin{align*}
\bigg(M(\na\times\mf F)\bigg)_{n}=\sum_{m,k,\ell}M_{nm}\vn_{m\ell k}\pl_{\ell}f_{k}.
\end{align*}
Finally, taking the curl operator on the vector $M(\na\times\mf F)$,
the $i$-th component of the resulted vector is 
\begin{align*}
\begin{aligned}\bigg(\na\times\big(M(\na\times\mf F)\big)\bigg)_{i} & =\sum_{j,n,m,k,\ell}\vn_{ijn}\pl_{j}\big(M_{nm}\vn_{m\ell k}\pl_{\ell}f_{k}\big)\\
 & =\sum_{j,n,m,k,\ell}\vn_{ijn}\vn_{m\ell k}\big((\pl_{j}M_{nm})\pl_{\ell}f_{k}+M_{nm}\pl_{j\ell}f_{k}\big)
\end{aligned}
\end{align*}
Thus 
\begin{align*}
\bigg(\na\times\big(M(\na\times\mf F)\big)\bigg)_{i}=\sum_{j,k,\ell}\wt C_{ijk\ell}\pl_{j\ell}f_{k}+\wt R_{i},
\end{align*}
where 
\begin{align*}
\wt C_{ijk\ell}:=\sum_{m,n}\vn_{ijn}\vn_{m\ell k}M_{nm},\mbox{\hspace{1cm}}\wt R_{i}:=\sum_{j,m,n,k,\ell}\vn_{ijn}\vn_{m\ell k}(\pl_{j}M_{nm})\pl_{\ell}f_{k}.
\end{align*}
Since 
\begin{align*}
\begin{aligned}\vn_{ijn}\vn_{m\ell k} & =\left|\begin{array}{ccc}
\da_{im} & \da_{i\ell} & \da_{ik}\\
\da_{jm} & \da_{j\ell} & \da_{jk}\\
\da_{nm} & \da_{n\ell} & \da_{nk}
\end{array}\right|\\
 & =\da_{im}\big(\da_{j\ell}\da_{nk}-\da_{n\ell}\da_{jk}\big)-\da_{i\ell}\big(\da_{jm}\da_{nk}-\da_{nm}\da_{jk}\big)+\da_{ik}\big(\da_{jm}\da_{n\ell}-\da_{nm}\da_{j\ell}\big),
\end{aligned}
\end{align*}
we can obtain 
\begin{align*}
\wt C_{ijk\ell} & =\sum_{mn}\bigg(\da_{im}\big(\da_{j\ell}\da_{nk}-\da_{n\ell}\da_{jk}\big)-\da_{i\ell}\big(\da_{jm}\da_{nk}-\da_{nm}\da_{jk}\big)\\
 & \qquad\qquad+\da_{ik}\big(\da_{jm}\da_{n\ell}-\da_{nm}\da_{j\ell}\big)\bigg)M_{nm}\\
 & =\big(\da_{j\ell}M_{ki}-\da_{jk}M_{\ell i}\big)-\da_{i\ell}M_{kj}+\da_{i\ell}\da_{jk}tr(M)+\da_{ik}M_{\ell j}-\da_{ik}\da_{j\ell}tr(M)\\
 & =\da_{j\ell}M_{ki}+\da_{ik}M_{\ell j}-\da_{jk}M_{\ell i}-\da_{i\ell}M_{kj}+\big(\da_{i\ell}\da_{jk}-\da_{ik}\da_{j\ell}\big)tr(M).
\end{align*}
\end{proof}
\begin{thm}
\label{Max-SE} Assume that $\mu$ is a smooth, positive scalar function
and $\en$ is a symmetric, positive definite matrix-valued function
with smooth entries. The eigenvalues of $\en$ are denoted by $\la_{1}(x),\la_{2}(x)$
and $\la_{3}(x)$. Assume there exist positive constants $\mu_{0}$,
$\Lambda,\la$ such that for all $x\in\Omega$ 
\begin{align*}
0<\mu(x)\le\mu_{0}
\end{align*}
\begin{align}
0<\la\leq\la_{1}(x)\leq\la_{2}(x)\leq\la_{3}(x)\leq\Lambda.\label{la}
\end{align}
Then (\ref{mA}) and (\ref{mB}) are uniformly strongly elliptic by
letting $M^{A}=m\mu^{-1}I$ and $M^{B}=m\mu^{-1}\en$, for arbitrary
positive constant $m$. Here $I$ denotes the $3\times3$ identity
matrix.\end{thm}
\begin{proof}
To see whether (\ref{mA}) and (\ref{mB}) are strongly elliptic,
we only have to check the leading order terms of (\ref{mA}) and (\ref{mB}).
We divide this proof into two parts, Part A and Part B, to deal with
the equation (\ref{mA}) for $A$ and the equation (\ref{mB}) for
$B$ respectively.
\begin{description}
\item [{\uline{Part\,A.}}] By Lemma \ref{id-curl}, 
\begin{align*}
\begin{aligned}\bigg(\mu\na & tr\big(M^{A}\na A\big)-\na\times\big(\ga^{-1}(\na\times A)\big)\bigg)_{i}\\
 & =\sum_{jk\ell}\mu\da_{ij}\pl_{j}\big(M_{\ell k}^{A}\pl_{\ell}A_{k}\big)-\sum_{jk\ell}\wt C_{ijk\ell}^{A}\pl_{j\ell}A_{k}-\wt R_{i}^{A}\\
 & =\sum_{jk\ell}\big(\mu\da_{ij}M_{\ell k}^{A}-\wt C_{ijk\ell}^{A}\big)\pl_{j\ell}A_{k}+\sum_{jk\ell}\mu\da_{ij}(\pl_{j}M_{\ell k}^{A})\pl_{\ell}A_{k}-\wt R_{i}^{A}\\
 & =\sum_{jk\ell}C_{ijk\ell}^{A}\pl_{j\ell}A_{k}+\sum_{jk\ell}\mu\da_{ij}(\pl_{j}M_{\ell k}^{A})\pl_{\ell}A_{k}-\wt R_{i}^{A},
\end{aligned}
\end{align*}
where $C_{ijk\ell}^{A}=\mu\da_{ij}M_{\ell k}^{A}-\wt C_{ijk\ell}^{A}$
are the coefficients of the leading order terms of (\ref{mA}) and
\begin{align*}
\begin{aligned}\wt C_{ijk\ell}^{A}=\da_{j\ell}(\en^{-1})_{ki}+\da_{ik}(\en^{-1})_{\ell j}-\da_{jk}(\en^{-1})_{\ell i}-\da_{i\ell}(\en^{-1})_{kj}+\big(\da_{i\ell}\da_{jk}-\da_{ik}\da_{j\ell}\big)tr(\en^{-1}).\end{aligned}
\end{align*}
Recall that (\ref{mA}) is called uniformly strongly elliptic in some
domain $\Omega$ if there exists a positive $c_{0}>0$ independent
of $x\in\Oa$ such that 
\begin{align}
\sum_{ijk\ell}C_{ijk\ell}^{A}(x)a_{i}a_{k}b_{j}b_{\ell}\geq c_{0}|\mf a|^{2}|\mf b|^{2}\label{useA}
\end{align}
for any $\mf a,\mf b\in\mb R^{3}$ and for all $x\in\Omega$. Now
\begin{align*}
\begin{aligned}\sum_{ijk\ell}C_{ijk\ell}^{A}a_{i}a_{k}b_{j}b_{\ell} & =\sum_{ijk\ell}\big(\mu\da_{ij}M_{\ell k}^{A}-\wt C_{ijk\ell}^{A}\big)a_{i}a_{k}b_{j}b_{\ell}\\
 & =\mu(\mf a\cdot\mf b)(\mf b^{T}M^{A}\mf a)\\
 & \quad-\sum_{ijk\ell}\bigg(\da_{j\ell}(\en^{-1})_{ki}+\da_{ik}(\en^{-1})_{\ell j}-\da_{jk}(\en^{-1})_{\ell i}\\
 & \qquad\qquad-\da_{i\ell}(\en^{-1})_{kj}+\big(\da_{i\ell}\da_{jk}-\da_{ik}\da_{j\ell}\big)tr(\en^{-1})\bigg)a_{i}a_{k}b_{j}b_{\ell}\\
 & =\mu(\mf a\cdot\mf b)(\mf b^{T}M^{A}\mf a)\\
 & \quad-\bigg(|\mf b|^{2}(\mf a^{T}\en^{-1}\mf a)+|\mf a|^{2}(\mf b^{T}\en^{-1}\mf b)-(\mf a\cdot\mf b)(\mf b^{T}\en^{-1}\mf a)\\
 & \qquad\qquad-(\mf a\cdot\mf b)(\mf a^{T}\en^{-1}\mf b)+tr(\en^{-1})(\mf a\cdot\mf b)^{2}-tr(\en^{-1})|\mf a|^{2}|\mf b|^{2}\bigg)\\
 & =tr(\en^{-1})|\mf a|^{2}|\mf b|^{2}-|\mf a|^{2}(\mf b^{T}\en^{-1}\mf b)-|\mf b|^{2}(\mf a^{T}\en^{-1}\mf a)-tr(\en^{-1})(\mf a\cdot\mf b)^{2}\\
 & \quad+2(\mf a\cdot\mf b)\big(\mf b^{T}\en^{-1}\mf a\big)+\mu(\mf a\cdot\mf b)\big(\mf b^{T}M^{A}\mf a\big)
\end{aligned}
\end{align*}
since $\en$ (and hence $\en^{-1}$) is symmetric. Let $S$ be the
orthogonal matrix such that $\en=S^{T}DS$, where $D=\mbox{diag}(\la_{1},\la_{2},\la_{3})$.
Thus $\en^{-1}=S^{T}D^{-1}S$. Also let $M^{A}=S^{T}N^{A}S$. By letting
$v=S\mf a/|\mf a|$ and $w=S\mf b/|\mf b|$, it's easy to see that
(\ref{useA}) holds for all $\mf a,\mf b\in\mb R^{3}$ iff 
\begin{align*}
 & tr(\en^{-1})-(\mf w^{T}D^{-1}\mf w)-(v^{T}D^{-1}\mf v)-tr(\en^{-1})(\mf v\cdot\mf w)^{2}\\
 & \quad+2(\mf v\cdot\mf w)\big(\mf w^{T}D^{-1}\mf v\big)+\mu(\mf v\cdot\mf w)\big(\mf w^{T}N^{A}\mf v\big)\ge c_{0}
\end{align*}
for all $\mf v,\mf w\in\mb R^{3}$ such that $|\mf v|=|\mf w|=1$.
Note that $tr(\en^{-1})=tr(D^{-1})=\la_{1}^{-1}+\la_{2}^{-1}+\la_{3}^{-1}$.
In summary, we find that (\ref{mA}) is uniformly strongly elliptic
on $\Oa$ iff 
\begin{align}
\inf_{\mf x\in\Oa}\bigg(\min_{|\mf v|=|\mf w|=1}F(\mf v,\mf w)\bigg)>0,\label{FA}
\end{align}
where 
\begin{align*}
F(\mf v,\mf w) & =\bigg(tr(D^{-1})-(\mf w^{T}D^{-1}\mf w)-(\mf v^{T}D^{-1}\mf v)-tr(D^{-1})(\mf v\cdot\mf w)^{2}\\
 & \quad+2(\mf v\cdot\mf w)\big(\mf w^{T}D^{-1}\mf v\big)\bigg)+\mu(\mf v\cdot\mf w)\big(\mf w^{T}N^{A}\mf v\big)\\
 & =:G(\mf v,\mf w)+\mu(\mf v\cdot\mf w)\big(\mf w^{T}N^{A}\mf v\big).
\end{align*}
We will show that 
\begin{align}
G(\mf v,\mf w)\ge\la_{3}^{-1}\big(1-(\mf v\cdot\mf w)^{2}\big)\label{GAlb}
\end{align}
under the constraints $|\mf v|=|\mf w|=1$. Then, by choosing $M^{A}=m\mu^{-1}I$
for some positive constant $m$, we also have $N^{A}=m\mu^{-1}I$,
and 
\begin{align*}
F(\mf v,\mf w) & =G(\mf v,\mf w)+m(\mf v\cdot\mf w)^{2}\\
 & \ge\la_{3}^{-1}\big(1-(\mf v\cdot\mf w)^{2}\big)+m(\mf v\cdot\mf w)^{2}\\
 & =\la_{3}^{-1}+(m-\la_{3}^{-1})(\mf v\cdot\mf w)^{2}.
\end{align*}
Now since $0\le(\mf v\cdot\mf w)^{2}\le1$, if $m\ge\la_{3}^{-1}$,
we have $F(\mf v,\mf w)\ge\la_{3}^{-1}$, while if $m<\la_{3}^{-1}$,
we have $F(\mf v,\mf w)\ge\la_{3}^{-1}+(m-\la_{3}^{-1})=m$. Remember
that $\la_{3}^{-1}(x)\ge\La^{-1}$ on $\Oa$, we conclude that $F(\mf v,\mf w)\ge\min(\La^{-1},m)$
for all $|\mf v|=|\mf w|=1$ and all $x\in\Oa$.\\
\\
It remains to show (\ref{GAlb}). For this, note that 
\begin{align*}
G(\mf v,\mf w)=\sum_{j=1,2,3}\la_{j}^{-1}\bigg(1-w_{j}^{2}-v_{j}^{2}-(\mf v\cdot\mf w)^{2}+2(\mf v\cdot\mf w)v_{j}w_{j}\bigg)=:\sum_{j}\la_{j}^{-1}K_{j}.
\end{align*}
We can prove $K_{j}\ge0$ as follows: Since $(\mf v\cdot\mf w)-v_{1}w_{1}=v_{2}w_{2}+v_{3}w_{3}$,
by Schwarz inequality we have 
\begin{align*}
|(\mf v\cdot\mf w)-v_{1}w_{1}|\le\sqrt{v_{2}^{2}+v_{3}^{2}}\sqrt{w_{2}^{2}+w_{3}^{2}}=\sqrt{1-v_{1}^{2}}\sqrt{1-w_{1}^{2}}.
\end{align*}
Taking square, we obtain 
\begin{align*}
(\mf v\cdot\mf w)^{2}-2(\mf v\cdot\mf w)v_{1}w_{1}+v_{1}^{2}w_{1}^{2}\le1-v_{1}^{2}-w_{1}^{2}+v_{1}^{2}w_{1}^{2},
\end{align*}
which means $K_{1}\ge0$. Similarly $K_{2},K_{3}\ge0$. As a consequence,
since $\la_{1}^{-1}\ge\la_{2}^{-1}\ge\la_{3}^{-1}$, we have 
\begin{align*}
G(\mf v,\mf w)\ge\la_{3}^{-1}(K_{1}+K_{2}+K_{3})=\la_{3}^{-1}\big(1-(\mf v\cdot\mf w)^{2}\big),
\end{align*}
which completes the proof of Part A.
\item [{\uline{Part\,B}.}] For (\ref{mB}), we have 
\begin{align}
\begin{aligned}\bigg(\ga\na & tr\big(M^{B}\na\mf B\big)-\na\times\big(\mu^{-1}(\na\times\mf B)\big)\bigg)_{i}\\
 & =\sum_{jk\ell}\ga_{ij}\pl_{j}\big(M_{\ell k}^{B}\pl_{\ell}B_{k}\big)-\sum_{jk\ell}\wt C_{ijk\ell}^{B}\pl_{j\ell}B_{k}-\wt R_{i}^{B}\\
 & =\sum_{jk\ell}\big(\ga_{ij}M_{\ell k}^{B}-\wt C_{ijk\ell}^{B}\big)\pl_{j\ell}B_{k}+\sum_{jk\ell}\ga_{ij}(\pl_{j}M_{\ell k}^{B})\pl_{\ell}B_{k}-\wt R_{i}^{B},
\end{aligned}
\label{mBprin}
\end{align}
where 
\begin{align*}
\begin{aligned}\wt C_{ijk\ell}^{B} & =\da_{j\ell}\mu^{-1}\da_{ki}+\da_{ik}\mu^{-1}\da_{\ell j}-\da_{jk}\mu^{-1}\da_{\ell i}\\
 & \qquad-\da_{i\ell}\mu^{-1}\da_{kj}+\big(\da_{i\ell}\da_{jk}-\da_{ik}\da_{j\ell}\big)tr(\mu^{-1}I)\\
 & =\mu^{-1}\Big(\da_{i\ell}\da_{jk}-\da_{ik}\da_{j\ell}\Big).
\end{aligned}
\end{align*}
Denote the coefficients of the leading order terms of (\ref{mBprin})
by $C_{ijk\ell}^{B}$, we have 
\begin{align*}
\begin{aligned}C_{ijk\ell}^{B}=\en_{ij}M_{\ell k}^{B}-\wt C_{ijk\ell}^{B}=\en_{ij}M_{\ell k}^{B}-\mu^{-1}\bigg(\da_{i\ell}\da_{jk}-\da_{ik}\da_{j\ell}\bigg).\end{aligned}
\end{align*}
By choosing $M^{B}=m\mu^{-1}\en$ we obtain 
\begin{align*}
\sum_{ijk\ell}C_{ijk\ell}^{B}a_{i}a_{k}b_{j}b_{\ell}=\mu^{-1}\bigg(m(\mf a^{T}\ga\mf b)^{2}-\Big((\mf a\cdot\mf b)^{2}-|\mf a|^{2}|\mf b|^{2}\Big)\bigg)
\end{align*}
for all $\mf a,\mf b\in\mb R^{3}$. Remember that $\en=S^{T}DS$.
Since we have assumed $\mu^{-1}\ge\mu_{0}$ for some positive constant
$\mu_{0}$, by letting $\mf v=S\mf a/|\mf a|$ and $\mf w=S\mf b/|\mf b|$
for $\mf a,\mf b\ne0$, we see to prove $C_{ijk\ell}^{B}a_{i}a_{k}b_{j}b_{\ell}\ge c_{0}|\mf a|^{2}|\mf b|^{2}$
for some constant $c_{0}>0$ is equivalent to prove 
\begin{align}
\inf_{\mf x\in\Oa}\min_{|\mf v|=|\mf w|=1}H(\mf v,\mf w)>0,\label{FB}
\end{align}
where $H(\mf v,\mf w)=m(\mf v^{T}D\mf w)^{2}+\big(1-(\mf v\cdot\mf w)^{2}\big)$.
Although (\ref{FB}) looks simpler than (\ref{FA}), we fail to find
a simple method as before to get a clear lower bound. Nevertheless,
it is also easy to see that (\ref{FB}) is true by continuity, as
follows: If $(\mf v\cdot\mf w)^{2}=1$, then $\mf v=\pm\mf w$, and
\begin{align*}
m(\mf v^{T}D\mf w)^{2}=m(\la_{1}v_{1}^{2}+\la_{2}v_{2}^{2}+\la_{3}v_{3}^{2})^{2}\ge m\la_{1}^{2}.
\end{align*}
By continuity, there exists $\vn>0$ such that for $0\le1-(\mf v\cdot\mf w)^{2}\le\vn$
we have $m(\mf v^{T}D\mf w)^{2}\ge m\la_{1}^{2}/2$. Thus for $0\le1-(\mf v\cdot\mf w)^{2}\le\vn$
we have $H(\mf v,\mf w)\ge m\la_{1}^{2}/2$. While for $1-(\mf v\cdot\mf w)^{2}>\vn$,
$H(\mf v,\mf w)>\vn$. Thus under the constraints $|\mf v|=|\mf w|=1$
we obtain 
\begin{align*}
H(\mf v,\mf w)\ge\min(m\la_{1}^{2}/2,\vn)\ge\min(m\la^{2}/2,\vn),
\end{align*}
where recall that $\la$ is the lower bound of $\la_{1}(x)$ on $\Oa$.
This completes the proof of Part B. 
\end{description}
\end{proof}
\begin{rem}
One can check that the $\wt C^{A}$ and $\wt C^{B}$ satisfy $\wt C_{ijk\ell}^{A}=\wt C_{k\ell ij}^{A}$
and $\wt C_{ijk\ell}^{B}=\wt C_{k\ell ij}^{B}$. And, by choosing
$M^{A}=m\mu^{-1}I$ and $M^{B}=m\mu^{-1}\en$ as above, the $C^{A}$
and $C^{B}$ also satisfy such symmetry. This additional property
is useful in the next section.
\end{rem}

\section{Construction of oscillating-decaying solutions}

In this section, we will use the reduction results in section 2 to
construct oscillating-decaying solutions of (\ref{ME}). From now
on, we suppose that $\mu>0$ is a $C^{\infty}$ scalar function and
$\epsilon$ is a $3\times3$ real positive definite matrix-valued
smooth functions (i.e. every entry is a real $C^{\infty}$ function)
and $E$ , $H$ satisfy 
\[
\begin{cases}
\nabla\times E-ik\mu H=0 & \mbox{ in }\Omega,\\
\nabla\times H+ik\epsilon E=0 & \mbox{ in }\Omega.
\end{cases}
\]

In order to obtain the oscillating-decaying solutions of $E$ and
$H$, we have to construct the oscillating-decaying solutions for
$A$ and $B$. We follow the proof in \cite{NUW2005(ODS)} to construct
the oscillating-decaying solutions for $A$ and $B$, but here we
need to derive higher derivatives for $A$ and $B$.

From \cite{NUW2005(ODS)}, we borrow several notations as follows.
Assume that $\Omega\subset\mathbb{R}^{3}$ is an open set with smooth
boundary and $\omega\in S^{2}$ is given. Let $\eta\in S^{2}$ and
$\zeta\in S^{2}$ be chosen so that $\{\eta,\zeta,\omega\}$ forms
an orthonormal system of $\mathbb{R}^{3}$. We then denote $x'=(x\cdot\eta,x\cdot\zeta)$.
Let $t\in\mathbb{R}$, $\Omega_{t}(\omega)=\Omega\cap\{x\cdot\omega>t\}$
and $\Sigma_{t}(\omega)=\Omega\cap\{x\cdot\omega=t\}$ be a non-empty
open set.
\begin{thm}
Given $\{\eta,\zeta,\omega\}$ an orthonormal system of $\mathbb{R}^{3}$,
$x'=(x\cdot\eta,x\cdot\zeta)$ and $t\in\mathbb{R}$. We set $\Omega_{t}(\omega)=\Omega\cap\{x\cdot\omega>t\}$
and $\Sigma_{t}(\omega)=\Omega\cap\{x\cdot\omega=t\}$, then We can
construct two types OD solutions for the Maxwell system in $\Omega_{t}(\omega)$
which can be useful for penetrable and impenetrable obstacles respectively.
There exist two solutions of (\ref{eq:1.4}) of the forms. The first
one is

\begin{equation}
\begin{cases}
E=F_{A}^{1}(x)e^{i\tau x\cdot\xi}e^{-\tau(x\cdot\omega-t)A_{t}^{A}(x')}b+\Gamma_{\chi_{t},b,t,N,\omega}^{A,1}(x,\tau)+r_{\chi_{t},b,t,N,\omega}^{A,1}(x,\tau) & \mbox{ in }\Omega_{t}(\omega),\\
H=F_{A}^{2}(x)e^{i\tau x\cdot\xi}e^{-\tau(x\cdot\omega-t)A_{t}^{A}(x')}b+\Gamma_{\chi_{t},b,t,N,\omega}^{A,2}(x,\tau)+r_{\chi_{t},b,t,N,\omega}^{A,2}(x,\tau) & \mbox{ in }\Omega_{t}(\omega),
\end{cases}\label{eq:3.1}
\end{equation}
where $F_{A}^{1}(x)=O(\tau)$, $F_{A}^{2}(x)=O(\tau^{2})$ are some
smooth functions and for $|\alpha|=j,$ $j=1,2$, we have 
\begin{equation}
\begin{cases}
\|\Gamma_{\chi_{t},b,t,N,\omega}^{A,j}(x,\tau)\|_{L^{2}(\Omega_{t}(\omega))}\leq c\tau^{|\alpha|-3/2}e^{-\tau(s-t)a_{A}},\\
\|r_{\chi_{t},b,t,N,\omega}^{A,j}(x,\tau)\|_{L^{2}(\Omega_{t}(\omega))}\leq c\tau^{j-N+1/2},
\end{cases}\label{eq:3.2}
\end{equation}
for some positive constants $a_{A}$ and $c$. The second one has
the form 
\begin{equation}
\begin{cases}
E=G_{B}^{2}(x)e^{i\tau x\cdot\xi}e^{-\tau(x\cdot\omega-t)A_{t}^{B}(x')}b+\Gamma_{\chi_{t},b,t,N,\omega}^{B,2}(x,\tau)+r_{\chi_{t},b,t,N,\omega}^{,B,2}(x,\tau) & \mbox{ in }\Omega_{t}(\omega),\\
H=G_{B}^{1}(x)e^{i\tau x\cdot\xi}e^{-\tau(x\cdot\omega-t)A_{t}^{B}(x')}b+\Gamma_{\chi_{t},b,t,N,\omega}^{B,1}(x,\tau)+r_{\chi_{t},b,t,N,\omega}^{B,1}(x,\tau) & \mbox{ in }\Omega_{t}(\omega),
\end{cases}\label{eq:3.3}
\end{equation}
where $G_{B}^{1}(x)=O(\tau)$,$G_{B}^{2}(x)=O(\tau^{2})$ are some
smooth functions and for $|\alpha|=j,$ $j=1,2$, we have 
\begin{equation}
\begin{cases}
\|\Gamma_{\chi_{t},b,t,N,\omega}^{B,j}(x,\tau)\|_{L^{2}(\Omega_{t}(\omega))}\leq c\tau^{|\alpha|-3/2}e^{-\tau(s-t)a_{B}},\\
\|r_{\chi_{t},b,t,N,\omega}^{B,j}(x,\tau)\|_{L^{2}(\Omega_{t}(\omega))}\leq c\tau^{j-N+1/2},
\end{cases}\label{eq:3.4}
\end{equation}
for some positive constants $a_{B}$ and $c$.\end{thm}
\begin{proof}
We want to find special solutions $A,B\in(C^{\infty}(\overline{\Omega_{t}(\omega)}\backslash\partial\Sigma_{t}(\omega))\cap C^{0}(\overline{\Omega_{t}(\omega)}))^{3}$
with $\tau\gg1$ satisfying Dirichlet boundary problems 
\begin{equation}
\begin{cases}
L_{A}A:=\mu\nabla tr(M^{A}\nabla A)-\nabla\times(\en^{-1}(\nabla\times A))+k^{2}\mu A=0 & \mbox{ in }\Omega_{t}(\omega)\\
A=e^{i\tau x\cdot\xi}\left\{ \chi_{t}(x')Q_{t}(x')b+\beta_{\chi_{t},t,b,N,\omega}^{A}\right\}  & \mbox{ on }\Sigma_{t}(\omega),
\end{cases}\label{eq:1.4}
\end{equation}
and 
\begin{equation}
\begin{cases}
L_{B}B:=\en\nabla tr(M^{B}\nabla B)-\nabla\times(\mu^{-1}(\nabla\times B))+k^{2}\en B=0 & \mbox{ in }\Omega_{t}(\omega)\\
B=e^{i\tau x\cdot\xi}\left\{ \chi_{t}(x')Q_{t}(x')b+\beta_{\chi_{t},t,b,N,\omega}^{B}\right\}  & \mbox{ on }\Sigma_{t}(\omega),
\end{cases}\label{eq:.1.5}
\end{equation}
where $\xi\in S^{2}$ lying in the span of $\{\eta,\zeta\}$ is chosen
and fixed, $\chi_{t}(x')\in C_{0}^{\infty}(\mathbb{R}^{2})$ with
supp$(\chi_{t})\subset\Sigma_{t}(\omega)$, $Q_{t}(x')$ is a nonzero
smooth function and $0\neq b\in\mathbb{C}^{3}$ and $N$ is some large
nature number. Moreover, $\beta_{\chi_{t},b,t,N,\omega}^{A}(x',\tau),\beta_{\chi_{t},b,t,N,\omega}^{B}(x',\tau)$
are smooth functions supported in supp($\chi_{t}$) satisfying: 
\[
\|\beta_{\chi_{t},b,t,N,\omega}^{A}(\cdot,\tau)\|_{L^{2}(\mathbb{R}^{2})}\leq c\tau^{-1},\mbox{ }\|\beta_{\chi_{t},b,t,N,\omega}^{B}(\cdot,\tau)\|_{L^{2}(\mathbb{R}^{2})}\leq c\tau^{-1}
\]
for some constant $c>0$. From now on, we use $c$ to denote a general
positive constant whose value may vary from line to line. As in \cite{NUW2005(ODS)},
$A,B$ satisfy second order strongly elliptic equations, then it can
be written as 
\[
\begin{cases}
A=A_{\chi_{t},b,t,N,\omega}=w_{\chi_{t},b,t,N,\omega}^{A}+r_{\chi_{t},b,t,N,\omega}^{A}\\
B=B_{\chi_{t},b,t,N,\omega}=w_{\chi_{t},b,t,N,\omega}^{B}+r_{\chi_{t},b,t,N,\omega}^{B}
\end{cases}
\]
with 
\begin{equation}
\begin{cases}
w_{\chi_{t},b,t,N,\omega}^{A}=\chi_{t}(x')Q_{t}e^{i\tau x\cdot\xi}e^{-\tau(x\cdot\omega-t)A_{t}^{A}(x')}b+\Gamma_{\chi_{t},b,t,N,\omega}^{A}(x,\tau)\\
w_{\chi_{t},b,t,N,\omega}^{B}=\chi_{t}(x')Q_{t}e^{i\tau x\cdot\xi}e^{-\tau(x\cdot\omega-t)A_{t}^{B}(x')}b+\Gamma_{\chi_{t},b,t,N,\omega}^{B}(x,\tau)
\end{cases}\label{eq:2.3}
\end{equation}
and $r_{\chi_{t}b,t,N,\omega}^{A},r_{\chi_{t}b,t,N,\omega}^{B}$ satisfying
\begin{equation}
\|r_{\chi_{t},b,t,N,\omega}^{A}\|_{H^{k}(\Omega_{t}(\omega))}\leq c\tau^{k-N+1/2},\mbox{ }\|r_{\chi_{t},b,t,N,\omega}^{B}\|_{H^{k}(\Omega_{t}(\omega))}\leq c\tau^{k-N+1/2},\label{eq:1.6}
\end{equation}
where $A_{t}^{A}(\cdot),\mbox{ }A_{t}^{B}(\cdot)$ are smooth matrix
functions with its real part Re$A_{t}^{A}(x')>0$, Re$A_{t}^{B}(x')>0$
and $\Gamma_{\chi_{t},b,t,N,\omega}^{A},\mbox{ }\Gamma_{\chi_{t},b,t,N,\omega}^{B}$
are a smooth functions supported in supp($\chi_{t}$) satisfying 
\begin{equation}
\begin{cases}
\|\partial_{x}^{\alpha}\Gamma_{\chi_{t},b,t,N,\omega}^{A}\|_{L^{2}(\Omega_{s}(\omega))}\leq c\tau^{|\alpha|-3/2}e^{-\tau(s-t)a_{A}}\\
\|\partial_{x}^{\alpha}\Gamma_{\chi_{t},b,t,N,\omega}^{B}\|_{L^{2}(\Omega_{s}(\omega))}\leq c\tau^{|\alpha|-3/2}e^{-\tau(s-t)a_{B}}
\end{cases}\label{eq:1.7}
\end{equation}
for $|\alpha|\in\mathbb{N}\cup\{0\}$ and $s\geq t$, where $a_{A},a_{B}>0$
are some constants depending on $A_{t}^{A}(x')$ and $A_{t}^{B}(x')$
respectively. We give details of the construction of $A$ and $B$
with the estimates (\ref{eq:2.3}) and (\ref{eq:1.6}) in the appendix.

In Appendix 6.1, we derive the explicit representation of $A$ and
$B$. Recall that $E$ and $H$ are represented in terms of $A$ and
$B$ as follows 
\begin{equation}
\begin{cases}
E=-\dfrac{i}{k}\gamma^{-1}\nabla\times(\mu^{-1}(\nabla\times B))-\gamma^{-1}(\nabla\times A),\\
H=\dfrac{i}{k}\mu^{-1}\nabla\times(\gamma^{-1}(\nabla\times A))-\mu^{-1}(\nabla\times B).
\end{cases}\label{eq:ODS1}
\end{equation}
Now, we can show that $(E,H)$ satisfies (\ref{eq:3.1}), (\ref{eq:3.2})
and we will use this form to prove Theorem \ref{main-thm} for the
penetrable case. Similarly, we can show that $(E,H)$ satisfies (\ref{eq:3.3}),
(\ref{eq:3.4}) in order to prove Theorem \ref{main-thm} for the
impenetrable case. All we need to do is to differentiate $A$ and
$B$ term by term componentwisely. For the main terms of $A$ and
$B$, we can differentiate $\chi_{t}(x')Q_{t}e^{i\tau x\cdot\xi}e^{-\tau(x\cdot\omega-t)A_{t}^{A}(x')}b$
and $\chi_{t}(x')Q_{t}e^{i\tau x\cdot\xi}e^{-\tau(x\cdot\omega-t)A_{t}^{B}(x')}b$
directly and it is easy to see that 
\[
\begin{cases}
\nabla\times A=\tau\widetilde{F_{A}}(x)e^{i\tau x\cdot\xi}e^{-\tau(x\cdot\omega-t)A_{t}^{A}(x')}b+\nabla\times\Gamma_{\chi_{t},b,t,N,\omega}^{A}(x,\tau)+\nabla\times r_{\chi_{t},b,t,N,\omega}^{A},\\
\nabla\times B=\tau\widetilde{F_{B}}(x)e^{i\tau x\cdot\xi}e^{-\tau(x\cdot\omega-t)A_{t}^{B}(x')}b+\nabla\times\Gamma_{\chi_{t},b,t,N,\omega}^{B}(x,\tau)+\nabla\times r_{\chi_{t},b,t,N,\omega}^{B},
\end{cases}
\]
where $\widetilde{F_{A}}(x)$ and $\widetilde{F_{B}}(x)$ are smooth
matrix-valued functions and support in supp$(\chi_{t}(x'))$. For
the penetrable obstacle case, we choose $A=w_{\chi_{t},b,t,N,\omega}^{A}+r_{\chi_{t},b,t,N,\omega}^{A}$
to be the oscillating-decaying solution satisfies $L_{A}A=0$ and
$B\equiv0$ (also satisfies $L_{B}0=0$) in $\Omega_{t}(\omega)$,
then (\ref{eq:ODS1}) will become to 
\[
\begin{cases}
E=-\gamma^{-1}(\nabla\times A),\\
H=\dfrac{i}{k}\mu^{-1}\nabla\times(\gamma^{-1}(\nabla\times A)),
\end{cases}
\]
which means 
\[
\begin{cases}
E=F_{A}^{1}(x)e^{i\tau x\cdot\xi}e^{-\tau(x\cdot\omega-t)A_{t}^{A}(x')}b+\Gamma_{\chi_{t},b,t,N,\omega}^{A,1}(x,\tau)+r_{\chi_{t},b,t,N,\omega}^{A,1}(x,\tau),\\
H=F_{A}^{2}(x)e^{i\tau x\cdot\xi}e^{-\tau(x\cdot\omega-t)A_{t}^{A}(x')}b+\Gamma_{\chi_{t},b,t,N,\omega}^{A,2}(x,\tau)+r_{\chi_{t},b,t,N,\omega}^{A,2}(x,\tau),
\end{cases}
\]
where $F_{A}^{1}(x)$, $F_{A}^{2}(x)$ are smooth functions consisting
$\mu(x)$, $\epsilon(x)$, $Q_{t}(x')$, $A_{t}^{A}(x')$ and their
curls (it can be seen by directly calculation). Moreover, by suitable
choice of $b$ (for example, we can choose $b\neq0$ is not parallel
to $\xi$), we will get $F_{A}^{1}(x)=O(\tau)$ and $F_{A}^{2}(x)=O(\tau^{2})$.
Moreover, $\Gamma_{\chi_{t},b,t,N,\omega}^{A,1}$ and $\Gamma_{\chi_{t},b,t,N,\omega}^{A,2}$
satisfy (\ref{eq:1.7}) for $|\alpha|=1$ and $|\alpha|=2$, respectively,
$r_{\chi_{t},b,t,N,\omega}^{A,1}$ and $r_{\chi_{t},b,t,N,\omega}^{A,1}$
satisfy (\ref{eq:1.6}) for $k=1$ and $k=2$, respectively. Similarly,
for the impenetrable obstacle case, we choose $A=0$ and $B=w_{\chi_{t},b,t,N,\omega}^{B}+r_{\chi_{t},b,t,N,\omega}^{B}$
in $\Omega_{t}(\omega)$, then 
\[
\begin{cases}
E=G_{B}^{2}(x)e^{i\tau x\cdot\xi}e^{-\tau(x\cdot\omega-t)A_{t}^{B}(x')}b+\Gamma_{\chi_{t},b,t,N,\omega}^{B,2}(x,\tau)+r_{\chi_{t},b,t,N,\omega}^{,B,2}(x,\tau),\\
H=G_{B}^{1}(x)e^{i\tau x\cdot\xi}e^{-\tau(x\cdot\omega-t)A_{t}^{B}(x')}b+\Gamma_{\chi_{t},b,t,N,\omega}^{B,1}(x,\tau)+r_{\chi_{t},b,t,N,\omega}^{B,1}(x,\tau),
\end{cases}
\]
where $G_{B}^{1}(x)=O(\tau)$ and $G_{B}^{2}(x)=O(\tau^{2})$ and
$\Gamma_{\chi_{t},b,t,N,\omega}^{B,j}$ satisfies (\ref{eq:1.7})
for $|\alpha|=j$ and $r_{\chi_{t},b,t,N,\omega}^{B,j}$ satisfies
(\ref{eq:1.7}) for $k=j$.
\end{proof}

\section{Runge approximation property}

\textcolor{black}{In this section, we derive the Runge approximation
property for the following} anisotropic Maxwell equation 
\[
\begin{cases}
\nabla\times E-ik\mu H=0\\
\nabla\times H+ik\epsilon E=0
\end{cases}\mbox{ in }\Omega,
\]
where $\mu$ is a smooth scalar function defined on $\Omega$ and
$\epsilon$ is a $3\times3$ smooth positive definite matrix. Recall
that 
\[
\mu(x)\geq\mu_{0}>0\mbox{ and }\sum_{i.j=1}^{3}\epsilon_{ij}(x)\xi_{i}\xi_{j}\geq\epsilon_{0}|\xi|^{2}\mbox{ }\forall\xi\in\mathbb{R}^{3}.
\]
If we set $u=\left(\begin{array}{c}
H\\
E
\end{array}\right)$ and 
\begin{equation}
L:=i\left(\begin{array}{cc}
\epsilon^{-1} & 0\\
0 & \mu^{-1}I_{3}
\end{array}\right)\left(\begin{array}{cc}
0 & \nabla\times\\
-\nabla\times & 0
\end{array}\right)+kI_{6},\label{eq:4.1}
\end{equation}
then we have 
\begin{equation}
Lu=0,\label{eq:Runge 1}
\end{equation}
where $I_{j}$ means $j\times j$ identity matrix for $j=3,6$. 
\begin{thm}
Let $D$ and $\Omega$ be two open bounded domains with $C^{\infty}$
boundary in $\mathbb{R}^{3}$ with $D\Subset\Omega$. If $u\in(H(curl,D))^{2}$
satisfies 
\[
Lu=0\mbox{ in }D.
\]
 Given any compact subset $K\subset D$ and any $\epsilon>0$, there
exists $U\in(H(curl,\Omega))^{2}$ such that 
\[
LU=0\mbox{ in }\Omega,
\]
and $\|U-u\|_{H(curl,K)}<\epsilon$, where $\|f\|_{H(curl,\Omega)}=\left(\|f\|_{L^{2}(\Omega)}+\|curlf\|_{L^{2}(\Omega)}\right)$. \end{thm}
\begin{proof}
The proof is standard and it is based on weak unique continuation
property for the anisotropic Maxwell system $L$ in (\ref{eq:4.1})
and the Hahn-Banach theorem. The unique continuation property of the
system $L$ is proved in \cite{leis2013initial}. For more details,
how to derive the Runge approximation property from the weak unique
continuation, we refer readers to \cite{lax1956stability}.
\end{proof}

\section{Proof of Theorem 1.1}

In this section, we want to use the Runge approximation property and
the OD solutions to prove Theorem \ref{main-thm}. We define $B$
to be an open ball in $\mathbb{R}^{3}$ such that $\overline{\Omega}\subset B$.
Assume that $\widetilde{\Omega}\subset\mathbb{R}^{3}$ is an open
Lipschitz domain with $\overline{B}\subset\widetilde{\Omega}$. Recall
we have set $\omega\in S^{2}$ and $\{\eta,\zeta,\omega\}$ forms
an orthonormal basis of $\mathbb{R}^{3}$ and $t_{0}=\inf_{x\in D}x\cdot\omega=x_{0}\cdot\omega$,
where $x_{0}=x_{0}(\omega)\in\partial D$.

\subsection{Penetrable Case}

For the anisotropic Maxwell's equation 
\begin{align}
\begin{cases}
\begin{aligned}\na\times E=ik\mu H\\
\na\times H=-ik\epsilon E\\
\mbox{div}(\epsilon E)=0\\
\mbox{div}(\mu H)=0,
\end{aligned}
\end{cases}.\label{ME-1}
\end{align}
for any $t\leq t_{0}$ and $\eta>0$ small enough, in section 3, we
have constructed 
\[
\begin{cases}
E_{t-\eta}=F_{A}^{1}(x)e^{i\tau x\cdot\xi}e^{-\tau(x\cdot\omega-(t-\eta))A_{t}^{A}(x')}b+\Gamma_{\chi_{t},b,t-\eta,N,\omega}^{A,1}(x,\tau)+r_{\chi_{t},b,t-\eta,N,\omega}^{A,1}(x,\tau),\\
H_{t-\eta}=F_{A}^{2}(x)e^{i\tau x\cdot\xi}e^{-\tau(x\cdot\omega-(t-\eta))A_{t}^{A}(x')}b+\Gamma_{\chi_{t},b,t-\eta,N,\omega}^{A,2}(x,\tau)+r_{\chi_{t},b,t-\eta,N,\omega}^{A,2}(x,\tau),
\end{cases}
\]
to be the oscillating-decaying solutions satisfying (\ref{ME-1})
in $B_{t-\eta}(\omega)=B\cap\{x|x\cdot\omega>t-\eta\}$, where $F_{A}^{1}(x)=O(\tau)$
and $F_{A}^{2}(x)=O(\tau^{2})$. Moreover, $\Gamma_{\chi_{t},b,t-\eta,N,\omega}^{A,1}$
and $\Gamma_{\chi_{t},b,t-\eta,N,\omega}^{A,2}$ satisfy (\ref{eq:1.7})
for $|\alpha|=1$ and $|\alpha|=2$, respectively, $r_{\chi_{t},b,t-\eta,N,\omega}^{A,1}$
and $r_{\chi_{t},b,t-\eta,N,\omega}^{A,1}$ satisfy (\ref{eq:1.6})
for $k=1$ and $k=2$, respectively. Similarly, we have 
\[
\begin{cases}
E_{t}=F_{A}^{1}(x)e^{i\tau x\cdot\xi}e^{-\tau(x\cdot\omega-t)A_{t}^{A}(x')}b+\Gamma_{\chi_{t},b,t,N,\omega}^{A,1}(x,\tau)+r_{\chi_{t},b,t,N,\omega}^{A,1}(x,\tau),\\
H_{t}=F_{A}^{2}(x)e^{i\tau x\cdot\xi}e^{-\tau(x\cdot\omega-t)A_{t}^{A}(x')}b+\Gamma_{\chi_{t},b,t,N,\omega}^{A,2}(x,\tau)+r_{\chi_{t},b,t,N,\omega}^{A,2}(x,\tau),
\end{cases}
\]
so be the oscillating-decaying solutions satisfying (\ref{ME-1})
in $B_{t}(\omega)=B\cap\{x|x\cdot\omega>t\}$, where $\Gamma_{\chi_{t},b,t,N,\omega}^{A,1}$
and $\Gamma_{\chi_{t},b,t,N,\omega}^{A,2}$ satisfy (\ref{eq:1.7})
for $|\alpha|=1$ and $|\alpha|=2$, respectively, $r_{\chi_{t},b,t,N,\omega}^{A,1}$
and $r_{\chi_{t},b,t,N,\omega}^{A,1}$ satisfy (\ref{eq:1.6}) for
$k=1$ and $k=2$, respectively. In fact, from the construction the
oscillating-decaying solutions and the property of continuous dependence
on parameters in ordinary differential equations in section 3, it
is not hard to see that for any $\tau$, 
\[
\begin{cases}
E_{t-\eta}\to E_{t}\\
H_{t-\eta}\to H_{t}
\end{cases}
\]
in $H^{2}(B_{t}(\omega))$ as $\eta$ tends to 0.

Note that $\overline{\Omega_{t}(\omega)}\subset B_{t-\eta}(\omega)$
for all $t\leq t_{0}$. By using the Runge approximation property,
we can see that there exists a sequence of functions $(E_{\eta,\ell},H_{\eta,\ell})$,
$\ell=1,2,\cdots$, such that 
\[
\begin{cases}
E_{\eta,\ell}\to E_{t-\eta}\\
H_{\eta,\ell}\to H_{t-\eta}
\end{cases}\mbox{ in }H(curl,B_{t}(\omega)),
\]
as $\ell\to\infty$, where $(E_{\eta,\ell},H_{\eta,\ell})$ satisfy
(\ref{ME-1}) in $\widetilde{\Omega}$ for all $\eta>0,\ell\in\mathbb{N}$.
Recall that the indicator function $I_{\rho}(\tau,t)$ was defined
by the formula: 
\[
I_{\rho}(\tau,t):=\lim_{\eta\to0}\lim_{\ell\to\infty}I_{\rho}^{\epsilon,\ell}(\tau,t),
\]
where 
\[
I_{\rho}^{\eta,\ell}(\tau,t):=ik\tau\int_{\partial\Omega}(\nu\times H_{\eta,\ell})\cdot(\overline{(\Lambda_{D}-\Lambda_{\emptyset})(\nu\times H_{\eta,\ell})}\times\nu)dS.
\]

We prove the Theorem 1.1 for the penetrable obstacle case. For the
anisotropic penetrable obstacle problem 
\begin{equation}
\begin{cases}
\nabla\times E-ik\mu H=0 & \mbox{ in }\Omega,\\
\nabla\times H+ik\epsilon E=0 & \mbox{ in }\Omega,\\
\nu\times H=f & \mbox{ on }\partial\Omega,
\end{cases}\label{eq:Inclusion1}
\end{equation}
where $k$ is not an eigenvalue of (\ref{eq:Inclusion1}). Moreover,
we assume $\mu$ is a positive smooth scalar function, $\epsilon=\epsilon_{0}(x)-\chi_{D}\epsilon_{D}(x)$,
where $\gamma_{0}$ is symmetric positive definite smooth matrix,
$\epsilon_{D}(x)$ is a symmetric smooth matrix with det$\epsilon_{D}(x)\neq0$
$\forall x\in D$ and $\chi_{D}=\begin{cases}
1 & \mbox{ }x\in D\\
0 & \mbox{ otherwise}
\end{cases}$. Moreover, we need $\epsilon=\epsilon(x)$ is a positive definite
matrix satisfying the uniform elliptic condition. Recall that when
$\epsilon(x)=\epsilon_{0}(x)$, we have constructed $E_{t}$ and $H_{t}$
which are oscillating-decaying solutions defined on the half space
for the anisotropic Maxwell's equation 
\begin{equation}
\begin{cases}
\nabla\times E-ik\mu H=0 & \mbox{ in }\Omega,\\
\nabla\times H+ik\epsilon E=0 & \mbox{ in }\Omega,
\end{cases}\label{eq:inclustion 1-1}
\end{equation}
and $\{(E_{\eta,\ell},H_{\eta,\ell})\}$ are sequence of functions
satisfying (\ref{eq:inclustion 1-1}) defined on the whole $\Omega$.
Therefore, we can define the boundary data $f_{\eta,\ell}=\nu\times H_{\eta,\ell}$
on $\partial\Omega$ and solve $(E,H)$ satisfies (\ref{eq:Inclusion1}).
Let $\widetilde{H_{\eta,\ell}}=H-H_{\eta,\ell}$ be the reflected
solution, then $\widetilde{H_{\eta,\ell}}$ satisfies 
\begin{equation}
\begin{cases}
\nabla\times(\epsilon^{-1}\nabla\times\widetilde{H_{\eta,\ell}})-k^{2}\mu\widetilde{H_{\eta,\ell}}=-\nabla\times((\epsilon^{-1}(x)-\epsilon_{0}^{-1}(x))\nabla\times H_{\eta,\ell})\mbox{ in }\Omega,\\
\nu\times\widetilde{H_{\eta,\ell}}=0\mbox{ on }\partial\Omega.
\end{cases}\label{eq:inclusion 2}
\end{equation}

\begin{lem}
We have the following estimates

1. 
\[
-\tau^{-1}I_{\rho}^{\eta,\ell}\geq\int_{D}[\epsilon(\epsilon^{-1}-\epsilon_{0}^{-1})^{-1}\epsilon_{0}^{-1}\nabla\times H_{\eta,\ell}]\cdot(\nabla\times\overline{H_{\eta,\ell}})dx-k^{2}\int_{\Omega}\mu|\widetilde{H_{\eta,\ell}}|^{2}dx.
\]

2. 
\[
\tau^{-1}I_{\rho}^{\eta,\ell}(\tau,t)\geq\int_{D}((\epsilon_{0}^{-1}-\epsilon^{-1})\nabla\times H_{\eta,\ell})\cdot(\nabla\times\overline{H_{\eta,\ell}})dx-k^{2}\int_{\Omega}\mu|\widetilde{H_{\eta,\ell}}|^{2}dx.
\]
\end{lem}
\begin{proof}
First, we need to prove the following identity 
\begin{eqnarray}
-\tau^{-1}I_{\rho}^{\eta,\ell}(\tau,t) & = & \int_{\Omega}\left((\epsilon^{-1}-\epsilon_{0}^{-1})\nabla\times H_{\eta,\ell}\right)\cdot(\nabla\times\overline{H_{\eta,\ell}})dx\nonumber \\
 &  & -\int_{\Omega}(\epsilon^{-1}\nabla\times\widetilde{H_{\eta,\ell}})\cdot(\nabla\times\overline{\widetilde{H_{\eta,\ell}}})dx-k^{2}\int_{\Omega}\mu|\widetilde{H_{\eta,\ell}}|^{2}dx.\label{iden1}
\end{eqnarray}
Multiplying $\overline{\widetilde{H_{\eta,l}}}$ in the equation (\ref{eq:inclusion 2})
and integrating by parts we have 
\begin{align*}
 & \int_{\Omega}(\epsilon^{-1}\nabla\times\widetilde{H_{\eta,\ell}})\cdot(\nabla\times\overline{\widetilde{H_{\eta,\ell}}})dx-k^{2}\int_{\Omega}\mu|\widetilde{H_{\eta,\ell}}|^{2}dx\\
 & +\int_{\Omega}((\epsilon^{-1}-\epsilon_{0}^{-1})\nabla\times H_{\eta,\ell})\cdot(\nabla\times\overline{\widetilde{H_{\eta,\ell}}})dx=0,
\end{align*}
\begin{align}
 & \int_{\Omega}(\epsilon^{-1}\nabla\times\widetilde{H_{\eta,\ell}})\cdot(\nabla\times\overline{\widetilde{H_{\eta,\ell}}})dx-k^{2}\int_{\Omega}\mu|\widetilde{H_{\eta,\ell}}|^{2}dx\nonumber \\
 & -\int_{\Omega}(\epsilon^{-1}-\epsilon_{0}^{-1})\nabla\times H_{\eta,\ell})\cdot(\nabla\times\overline{H_{\eta,\ell}})dx\label{eq:lem5.3-1}\\
= & -\int_{\Omega}((\epsilon^{-1}-\epsilon_{0}^{-1})\nabla\times H_{\eta,\ell})\cdot(\nabla\times\overline{H})dx.
\end{align}
On the other hand, $H(x)$ satisfies 
\begin{equation}
\nabla\times(\epsilon^{-1}(x)\nabla\times H(x))-k^{2}\mu H(x)=0,\label{eq:inclusion 3}
\end{equation}
then multiply by $H_{\eta,l}(x)$ in the equation (\ref{eq:inclusion 3})
and integrating by parts we have
\begin{eqnarray}
\int_{\Omega}((\epsilon^{-1}-\epsilon_{0}^{-1})\nabla\times H_{\eta,\ell})\cdot(\nabla\times\overline{H})dx & = & \int_{\partial\Omega}(\epsilon^{-1}\nabla\times\overline{H})\cdot(\nu\times H_{\eta,\ell})ds\nonumber \\
 &  & -\int_{\partial\Omega}(\epsilon_{0}^{-1}\nabla\times H_{\eta,\ell})\cdot(\nu\times\overline{H})ds\label{eq:lem5.3-2}
\end{eqnarray}
Thus, combine (\ref{eq:lem5.3-1}), (\ref{eq:lem5.3-2}) and $\int_{\partial\Omega}(\nu\times\overline{H_{\eta,\ell}})\cdot(\epsilon_{0}^{-1}\nabla\times H_{\eta,\ell})ds$
is real, then we have 
\begin{align}
 & \int_{\Omega}(\epsilon^{-1}\nabla\times\widetilde{H_{\eta,\ell}})\cdot(\nabla\times\overline{\widetilde{H_{\eta,\ell}}})dx-k^{2}\int_{\Omega}\mu|\widetilde{H_{\eta,\ell}}|^{2}dx\nonumber \\
 & -\int_{\Omega}((\epsilon^{-1}-\epsilon_{0}^{-1})\nabla\times H_{\eta,\ell})\cdot(\nabla\times\overline{H_{\eta,\ell}})dx\\
= & \int_{\partial\Omega}(\nu\times H_{\eta,\ell})\cdot(\epsilon^{-1}\nabla\times\overline{H})ds-\int_{\partial\Omega}(\nu\times\overline{H})\cdot(\epsilon_{0}^{-1}\nabla\times H_{\eta,\ell})ds\nonumber \\
= & \int_{\partial\Omega}(\nu\times H_{\eta,\ell})\cdot(\epsilon^{-1}\nabla\times\overline{H})ds-\int_{\partial\Omega}(\nu\times\overline{H_{\eta,\ell}})\cdot(\epsilon_{0}^{-1}\nabla\times H_{\eta,\ell})ds\nonumber \\
= & \int_{\partial\Omega}(\nu\times H_{\eta,\ell})\cdot(\epsilon^{-1}\overline{\nabla\times H})ds-\int_{\partial\Omega}(\nu\times H_{\eta,\ell})\cdot(\epsilon_{0}^{-1}\overline{\nabla\times H_{\eta,\ell}})ds\nonumber \\
= & \int_{\partial\Omega}(\nu\times H_{\eta,\ell})\cdot[\overline{-ikE+ikE_{\eta,\ell}}]ds\nonumber \\
= & ik\int_{\partial\Omega}(\nu\times H_{\eta,\ell})\cdot[\overline{(\Lambda_{D}-\Lambda_{\emptyset})(\nu\times H_{\eta,\ell})}\times\nu]ds\nonumber \\
= & \tau^{-1}I_{\rho}^{\eta,\ell}.\label{eq:lem5.3-2-2}
\end{align}

Second, we show the following identity 
\begin{align}
 & \int_{\Omega}(\epsilon_{0}^{-1}\nabla\times\widetilde{H_{\eta,\ell}})\cdot(\nabla\times\overline{\widetilde{H_{\eta,\ell}}})dx-k^{2}\int_{\Omega}\mu|\widetilde{H_{\eta,\ell}}|^{2}dx\label{iden2}\\
 & +\int_{\Omega}((\epsilon^{-1}(x)-\epsilon_{0}^{-1}(x))\nabla\times H)\cdot(\nabla\times\overline{H})dx\nonumber \\
= & -\tau^{-1}I_{\rho}^{\eta,\ell}.\nonumber 
\end{align}
Replacing $H_{\eta,\ell}(x)$ by $H(x)-\widetilde{H_{\eta,\ell}}(x)$
in the equation (\ref{eq:inclusion 2}), then we have
\begin{equation}
\nabla\times\left((\epsilon^{-1}-\epsilon_{0}^{-1})\nabla\times H\right)+\nabla\times\left(\epsilon_{0}^{-1}\nabla\times\widetilde{H_{\eta,\ell}}\right)-k^{2}\mu\widetilde{H_{\eta,\ell}}=0\mbox{ in }\Omega.\label{eq:lem5.3-3}
\end{equation}
Multiplying $\overline{\widetilde{H_{\eta,l}}}(x)$ in the equation
(\ref{eq:lem5.3-3}) and using integration by parts we have

\begin{align}
 & \int_{\Omega}\left((\epsilon^{-1}-\epsilon_{0}^{-1})\nabla\times H\right)\cdot\left(\nabla\times\overline{\widetilde{H_{\eta,\ell}}}\right)dx\nonumber \\
 & +\int_{\Omega}\left(\epsilon_{0}^{-1}\nabla\times\widetilde{H_{\eta,\ell}}\right)\cdot\left(\nabla\times\overline{\widetilde{H_{\eta,\ell}}}\right)dx-k^{2}\int_{\Omega}\mu\left|\widetilde{H_{\eta,\ell}}\right|^{2}dx=0,\label{eq:lem5.3-4}
\end{align}
since $\nu\times\widetilde{H_{\eta,l}}=0$ on $\partial\Omega$. Then
we can write equation (\ref{eq:lem5.3-4}) to be 
\begin{align}
 & \int_{\Omega}\left(\epsilon_{0}^{-1}\nabla\times\widetilde{H_{\eta,\ell}}\right)\cdot\left(\nabla\times\overline{\widetilde{H_{\eta,\ell}}}\right)dx-k^{2}\int_{\Omega}\mu\left|\widetilde{H_{\eta,\ell}}\right|^{2}dx\nonumber \\
 & +\int_{\Omega}\left((\epsilon^{-1}-\epsilon_{0}^{-1})\nabla\times H\right)\cdot(\nabla\times\overline{H})dx\nonumber \\
= & \int_{\Omega}\left((\epsilon^{-1}-\epsilon_{0}^{-1})\nabla\times H\right)\cdot(\nabla\times\overline{H_{\eta,\ell}})dx.\label{eq:lem5.3-5}
\end{align}
Eliminating $H(x)$ by $\widetilde{H_{\eta,l}}(x)+H_{\eta,l}(x)$
in (\ref{eq:lem5.3-5}) we have 
\begin{align}
 & \int_{\Omega}\left(\epsilon_{0}^{-1}\nabla\times\widetilde{H_{\eta,\ell}}\right)\cdot\left(\nabla\times\overline{\widetilde{H_{\eta,\ell}}}\right)dx-k^{2}\int_{\Omega}\mu\left|\widetilde{H_{\eta,\ell}}\right|^{2}dx\nonumber \\
 & +\int_{\Omega}\left((\epsilon^{-1}-\epsilon_{0}^{-1})\nabla\times H\right)\cdot(\nabla\times\overline{H})dx\nonumber \\
= & \int_{\Omega}((\epsilon^{-1}(x)-\epsilon_{0}^{-1}(x))\nabla\times H_{\eta,\ell})\cdot(\nabla\times\overline{H_{\eta,\ell}})dx\nonumber \\
 & +\int_{\Omega}((\epsilon^{-1}(x)-\epsilon_{0}^{-1}(x))\nabla\times\widetilde{H_{\eta,\ell}})\cdot(\nabla\times\overline{H_{\eta,\ell}})dx\label{eq:lem5.3-6}
\end{align}
Again from (\ref{eq:inclusion 2}) \textcolor{black}{and by} taking
the complex conjugate, we can write 
\begin{equation}
\nabla\times(\epsilon^{-1}\nabla\times\overline{\widetilde{H_{\eta,\ell}}})-k^{2}\mu\overline{\widetilde{H_{\eta,\ell}}}+\nabla\times((\epsilon^{-1}(x)-\epsilon_{0}^{-1}(x))\nabla\times\overline{H_{\eta,\ell}})=0.\label{eq:lem5.3-7}
\end{equation}
Multiplying by $\widetilde{H_{\eta,l}}(x)$ in the equation (\ref{eq:lem5.3-7})
and using integration by parts we have 
\begin{eqnarray}
\int_{\Omega}(\epsilon^{-1}\nabla\times\overline{\widetilde{H_{\eta,\ell}}})\cdot(\nabla\times\widetilde{H_{\eta,\ell}})dx-k^{2}\int_{\Omega}\mu|\widetilde{H_{\eta,\ell}}|^{2}dx\nonumber \\
+\int_{\Omega}((\epsilon^{-1}(x)-\epsilon_{0}^{-1}(x))\nabla\times\overline{H_{\eta,\ell}})\cdot(\nabla\times\widetilde{H_{\eta,\ell}})dx=0.\label{eq:lem5.3-8}
\end{eqnarray}
Then from the equations (\ref{eq:lem5.3-6}), (\ref{eq:lem5.3-8})
\textcolor{black}{and the first identity (\ref{iden1}),} we can obtain
\begin{align}
 & \int_{\Omega}(\epsilon_{0}^{-1}\nabla\times\widetilde{H_{\eta,\ell}})\cdot(\nabla\times\overline{\widetilde{H_{\eta,\ell}}})dx-k^{2}\int_{\Omega}\mu|\widetilde{H_{\eta,\ell}}|^{2}dx\nonumber \\
 & +\int_{\Omega}((\epsilon^{-1}(x)-\epsilon_{0}^{-1}(x))\nabla\times H)\cdot(\nabla\times\overline{H})dx\nonumber \\
= & \int_{\Omega}((\epsilon^{-1}(x)-\epsilon_{0}^{-1}(x))\nabla\times H_{\eta,\ell})\cdot(\nabla\times\overline{H_{\eta,\ell}})dx\nonumber \\
 & -\int_{\Omega}(\epsilon^{-1}\nabla\times\overline{\widetilde{H_{\eta,\ell}}})\cdot(\nabla\times\widetilde{H_{\eta,\ell}})dx+k^{2}\int_{\Omega}\mu|\widetilde{H_{\eta,\ell}}|^{2}dx\nonumber \\
= & -\tau^{-1}I_{\rho}^{\eta,\ell}.\label{eq:lem5.3-9}
\end{align}
 Combine (\ref{eq:lem5.3-9}) with the formula
\begin{align*}
 & (\epsilon_{0}^{-1}\nabla\times\widetilde{H_{\eta,\ell}})\cdot(\nabla\times\overline{\widetilde{H_{\eta,\ell}}})+((\epsilon^{-1}-\epsilon_{0}^{-1})\nabla\times H)\cdot(\nabla\times\overline{H})\\
= & ((\epsilon^{-1}-\epsilon_{0}^{-1})\nabla\times H)\cdot\nabla\times\overline{H}+\epsilon_{0}^{-1}(\na\times H)\cdot(\na\times\ol H)\\
 & -2\mbox{Re}\left\{ \epsilon_{0}^{-1}\nabla\times H\cdot\nabla\times\overline{H_{\eta,\ell}}\right\} +\epsilon_{0}^{-1}\nabla\times H_{\eta,\ell}\cdot\nabla\times\overline{H_{\eta,\ell}}\\
= & \epsilon^{-1}(\na\times H)\cdot(\na\times\ol H)-2\mbox{Re}\left\{ \epsilon_{0}^{-1}\nabla\times H\cdot\nabla\times\overline{H_{\epsilon,l}}\right\} +\epsilon_{0}^{-1}\nabla\times H_{\epsilon,\ell}\cdot\nabla\times\overline{H_{\epsilon,\ell}}\\
= & \left[\epsilon^{-\fc 12}\na\times H-\epsilon^{\fc 12}\epsilon_{0}^{-1}\left(\na\times\overline{H_{\eta,\ell}}\right)\right]\cdot\left[\overline{\epsilon^{-\fc 12}\na\times H-\epsilon^{\fc 12}\epsilon_{0}^{-1}\left(\na\times\overline{H_{\eta,\ell}}\right)}\right]\\
 & -\left[\epsilon^{\fc 12}\epsilon_{0}^{-1}\left(\na\times\ol{H_{\eta,\ell}}\right)\right]\cdot\ol{\left[\epsilon^{\fc 12}\epsilon_{0}^{-1}\left(\na\times\ol{H_{\eta,\ell}}\right)\right]}+\epsilon_{0}^{-1}\nabla\times H_{\eta,\ell}\cdot\nabla\times\overline{H_{\eta,\ell}}\\
= & \left[\epsilon^{-\fc 12}\na\times H-\epsilon^{\fc 12}\epsilon_{0}^{-1}\left(\na\times\ol{H_{\eta,\ell}}\right)\right]\cdot\ol{\left[\epsilon^{-\fc 12}\na\times H-\epsilon^{\fc 12}\epsilon_{0}^{-1}\left(\na\times\ol{H_{\eta,\ell}}\right)\right]}\\
 & +\left(\epsilon_{0}^{-1}-\epsilon\epsilon_{0}^{-2}\right)\left(\na\times H_{\eta,\ell}\right)\cdot\left(\na\times\ol{H_{\eta,\ell}}\right)\\
\geq & [\left(I-\epsilon\epsilon_{0}^{-1}\right)\epsilon_{0}^{-1}\nabla\times H_{\eta,\ell}]\cdot(\nabla\times\overline{H_{\eta,\ell}})\\
\geq & [\epsilon(\epsilon^{-1}-\epsilon_{0}^{-1})^{-1}\epsilon_{0}^{-1}\nabla\times H_{\eta,\ell}]\cdot(\nabla\times\overline{H_{\eta,\ell}})
\end{align*}
and note that 
\[
\left[\epsilon^{-\fc 12}\na\times H-\epsilon^{\fc 12}\epsilon_{0}^{-1}\left(\na\times\ol{H_{\eta,\ell}}\right)\right]\cdot\ol{\left[\epsilon^{-\fc 12}\na\times H-\epsilon^{\fc 12}\epsilon_{0}^{-1}\left(\na\times\ol{H_{\eta,\ell}}\right)\right]}\geq0.
\]
Therefore, we get 
\[
-\tau^{-1}I_{\rho}^{\eta,\ell}\geq\int_{D}[\epsilon(\epsilon^{-1}-\epsilon_{0}^{-1})^{-1}\epsilon_{0}^{-1}\nabla\times H_{\eta,\ell}]\cdot(\nabla\times\overline{H_{\eta,\ell}})dx-k^{2}\int_{\Omega}\mu|\widetilde{H_{\eta,\ell}}|^{2}dx
\]
which finished the part 1 of lemma 4.1. Finally, again from (\ref{eq:lem5.3-2-2}),
we have 
\[
\tau^{-1}I_{\rho}^{\eta,\ell}\geq\int_{\Omega}((\epsilon_{0}^{-1}-\epsilon^{-1})\nabla\times H_{\eta,\ell})\cdot(\nabla\times\overline{H_{\eta,\ell}})dx-k^{2}\int_{\Omega}\mu|\widetilde{H_{\eta,\ell}}|^{2}dx.
\]
\end{proof}
\begin{rem}
The first inequality will be used when $\left(\epsilon^{-1}-\epsilon_{0}^{-1}\right)$
is strictly positive definite, i.e.
\[
\xi\cdot(\epsilon^{-1}-\epsilon_{0}^{-1})\xi\geq\Lambda|\xi|^{2}\mbox{ for all }\xi\in\mathbb{R}^{3}\mbox{ and for some }\Lambda>0;
\]
and the second inequality will be used when $\left(\epsilon_{0}^{-1}-\epsilon^{-1}\right)$
is strictly positive definite, i.e. 
\[
\xi\cdot(\epsilon_{0}^{-1}-\epsilon^{-1})\xi\geq\lambda|\xi|^{2}\mbox{ for all }\xi\in\mathbb{R}^{3}\mbox{ and for some }\la>0.
\]

\end{rem}
Now, our work is to estimate the lower order term $\widetilde{H_{\eta,\ell}}$.

\subsubsection{Estimate of the lower order term $\widetilde{H_{\eta,\ell}}$}
\begin{prop}
Assume $\Omega$ is a smooth domain and $D\Subset\Omega$. Then there
exist a positive constant $C$ and $\delta>0$ such that 
\[
\|\widetilde{H_{\eta,\ell}}\|_{L^{2}(\Omega)}\leq C\|\nabla\times H_{\eta,\ell}\|_{L^{p}(D)}
\]
for every $p\in(\max\{\dfrac{4}{3},\dfrac{2+\delta}{1+\delta}\},2]$.\end{prop}
\begin{proof}
We follow the proof of the proposition 3.2 in \cite{kar2014reconstruction}.
Fix $l\in\mathbb{N}$ and we set $f:=-(\epsilon^{-1}-\epsilon_{0}^{-1})(\nabla\times H_{\eta,\ell})$,
$g=0$. Note that, $\epsilon^{-1}-\epsilon_{0}^{-1}=\epsilon^{-1}(\epsilon_{D}\chi_{D})\epsilon_{0}^{-1}$
is supported in $D$. Then the reflected solution $\widetilde{H_{\eta,\ell}}$
satisfies
\begin{equation}
\begin{cases}
\nabla\times(\epsilon^{-1}\nabla\widetilde{H_{\eta,\ell}})-k^{2}\mu\widetilde{H_{\eta,\ell}}=-\nabla\times((\epsilon^{-1}(x)-\epsilon_{0}^{-1}(x))\nabla\times H_{\eta,\ell})\mbox{ in }\Omega,\\
\nu\times\widetilde{H_{\eta,\ell}}=0\mbox{ on }\partial\Omega.
\end{cases}\label{eq:prop4.1}
\end{equation}
From the $L^{p}$ estimate (Theorem 6.6), if we consider the following
problem
\[
\begin{cases}
\nabla\times(\epsilon^{-1}\nabla\times U)+\epsilon_{\max}^{-1}U=\nabla\times f & \mbox{ in }\Omega,\\
\nu\times U=0 & \mbox{ on }\partial\Omega,
\end{cases}
\]
has a unique solution in $H_{0}^{1,q}(curl,\Omega)$, where $\epsilon_{\max}^{-1}$
is the maximum value among all eigenvalues of the matrix $\epsilon^{-1}(x)$
in the region $\overline{\Omega}$. Moreover, we have the estimate
\begin{equation}
\|U\|_{L^{p}(\Omega)}+\|\nabla\times U\|_{L^{p}(\Omega)}\leq C\|f\|_{L^{p}(\Omega)}\label{eq:prop4.2}
\end{equation}
for $p\in(\dfrac{2+\delta}{1+\delta},2]$ for some $\delta>0$ which
depends only on $\Omega$. Now, we set $\Pi_{\eta,\ell}=\widetilde{H_{\eta,\ell}}-U$,
then $\Pi_{\eta,\ell}$ satisfies 
\begin{equation}
\begin{cases}
\nabla\times(\epsilon^{-1}\nabla\Pi_{\eta,\ell})-k^{2}\mu\Pi_{\eta,\ell}=(k^{2}\mu+\epsilon_{\max}^{-1})U\mbox{ in }\Omega,\\
\nu\times\Pi_{\eta,\ell}=0\mbox{ on }\partial\Omega. & \mbox{}
\end{cases}\label{eq:LOT1}
\end{equation}
By the well-posedness of (\ref{eq:LOT1}) in $H(curl,\Omega)$ for
the anisotropic Maxwell's equation (see Appendix), we have 
\begin{equation}
\|\Pi_{\eta,\ell}\|_{L^{2}(\Omega)}+\|\nabla\times\Pi_{\eta,\ell}\|_{L^{2}(\Omega)}\leq C\|U\|_{L^{2}(\Omega)}\label{eq:prop4.3}
\end{equation}
if $k$ is not an eigenvalue. Moreover, for $p\leq2$, it is to see
that 
\[
\|\Pi_{\eta,\ell}\|_{L^{p}(\Omega)}+\|\nabla\times\Pi_{\eta,\ell}\|_{L^{p}(\Omega)}\leq C\|U\|_{L^{2}(\Omega)}.
\]
Following the proof in the proposition 3.2 in \cite{KS2014} again,
we denote $B_{\frac{1}{p}}^{p,2}(\Omega)$ to be the Sobolev-Besov
space, then we have $U\in B_{\frac{1}{p}}^{p,2}(\Omega)$ and the
inclusion map $B_{\frac{1}{p}}^{p,2}(\Omega)\to L^{2}(\Omega)$ is
continuous for $p\in(\tfrac{4}{3},2]$. Moreover, since $\nabla\times U=0$
and $\nu\times U=0$ on $\partial\Omega$ and use Lemma 7.6 ( property
5 in the appendix of \cite{KS2014}), we have the estimate 
\begin{equation}
\|U\|_{L^{2}(\Omega)}\le C\|U\|_{B_{\frac{1}{p}}^{p,2}(\Omega)}\leq C\{\|U\|_{L^{p}(\Omega)}+\|\nabla\times U\|_{L^{p}(\Omega)}\}\label{eq:prop4.4}
\end{equation}
for $p\in(\tfrac{4}{3},2]$. Combining (\ref{eq:prop4.2}), (\ref{eq:prop4.3})
and (\ref{eq:prop4.4}), we obtain 
\begin{equation}
\|\Pi_{\eta,\ell}\|_{L^{p}(\Omega)}+\|\nabla\times\Pi_{\eta,\ell}\|_{L^{p}(\Omega)}\leq C\|f\|_{L^{p}(\Omega)}\label{eq:prop4.5}
\end{equation}
for $p\in(\max\{\tfrac{4}{3},\tfrac{2+\delta}{1+\delta}\},2]$. Since
$\widetilde{H_{\eta,\ell}}=\Pi_{\eta,\ell}+U$, by using (\ref{eq:prop4.2})
and (\ref{eq:prop4.5}), we have 
\begin{equation}
\|\widetilde{H_{\eta,\ell}}\|_{L^{p}(\Omega)}+\|\nabla\times\widetilde{H_{\eta,\ell}}\|_{L^{p}(\Omega)}\leq C\|f\|_{L^{p}(\Omega)}.\label{eq:prop4.6}
\end{equation}
Since $\nu\times\widetilde{H_{\eta,l}}=0$ on $\partial\Omega$, we
use the Lemma 7.6 again, then we can obtain 
\begin{eqnarray}
\|\widetilde{H_{\eta,\ell}}\|_{L^{2}(\Omega)} & \le & C\|\widetilde{H_{\eta,\ell}}\|_{B_{\frac{1}{p}}^{p,2}(\Omega)}\nonumber \\
 & \leq & C\{\|\widetilde{H_{\eta,\ell}}\|_{L^{p}(\Omega)}+\|\nabla\times\widetilde{H_{\eta,\ell}}\|_{L^{p}(\Omega)}+\|\nabla\cdot\widetilde{H_{\eta,\ell}}\|_{L^{p}(\Omega)}\}.\label{eq:prop4.7}
\end{eqnarray}
In addition, from (\ref{eq:prop4.1}), it is easy to see $0=\nabla\cdot(\mu\widetilde{H_{\eta,\ell}})=\nabla\mu\cdot\widetilde{H_{\eta,\ell}}+\mu(\nabla\cdot\widetilde{H_{\eta,\ell}})$,
then we have 
\begin{equation}
\|\nabla\cdot\widetilde{H_{\eta,\ell}}\|_{L^{p}(\Omega)}\leq\dfrac{\|\nabla\mu\|_{L^{\infty}(\Omega)}}{\|\mu\|_{L^{\infty}(\Omega)}}\|\widetilde{H_{\eta,\ell}}\|_{L^{p}(\Omega)}.\label{eq:prop4.8}
\end{equation}
Finally, use (\ref{eq:prop4.6}), (\ref{eq:prop4.7}) and (\ref{eq:prop4.8}),
we will get 
\begin{eqnarray}
\|\widetilde{H_{\eta,\ell}}\|_{L^{2}(\Omega)} & \leq & C\{\|\widetilde{H_{\eta,\ell}}\|_{L^{p}(\Omega)}+\|\nabla\times\widetilde{H_{\eta,\ell}}\|_{L^{p}(\Omega)}\}\nonumber \\
 & \leq & C\|f\|_{L^{p}(\Omega)}\nonumber \\
 & \leq & C\|\nabla\times H_{\eta,\ell}\|_{L^{p}(D)}.\label{eq:prop4.9}
\end{eqnarray}
\end{proof}
\begin{rem}
In the reconstruction scheme, we need to take $\limsup_{\ell\to\infty}$
for (\ref{eq:prop4.9}) on both sides and $H_{t+\eta}\to H_{t}$ in
$H(curl,\Omega_{t}(\omega))$ as $\eta\to0$, then we have 
\[
\lim_{\eta\to0}\limsup_{\ell\to\infty}\|\widetilde{H_{\eta,\ell}}\|_{L^{2}(\Omega)}\leq C\|\nabla\times H_{t}\|_{L^{p}(D)},
\]
for $p\in(\dfrac{4}{3},2]$. Moreover, if (\ref{eq:prop4.1}) is written
as the following form 
\[
\begin{cases}
\nabla\times(\en^{-1}\nabla\widetilde{H})-k^{2}\mu\widetilde{H}=-\nabla\times((\en^{-1}(x)-\en_{0}^{-1}(x))\nabla\times H_{0})+k^{2}(\en-\en_{0})H_{0}.\mbox{ in }\Omega,\\
\nu\times\widetilde{H}=0\mbox{ on }\partial\Omega,
\end{cases}
\]
and we can derive the following estimate by using the same method
in the proof of the Proposition 5.3, then the estimate (\ref{eq:prop4.9})
will be 
\[
\|\widetilde{H}\|_{L^{2}(\Omega)}\leq C\{\|\nabla\times H_{0}\|_{L^{p}(D)}+\|H_{0}\|_{L^{2}(\Omega)}\},
\]
for $p\in(\dfrac{4}{3},2]$.
\end{rem}
In view of the lower bound, we need to introduce the sets $D_{j,\delta}\subset D$,
$D_{\delta}\subset D$ in the following. Recall that $h_{D}(\rho)=\inf_{x\in D}x\cdot\rho$
and $t_{0}=h_{D}(\rho)=x_{0}\cdot\rho$ for some $x_{0}\in\partial D$.
$\forall\alpha\in\partial D\cap\{x\cdot\rho=h_{D}(\rho)\}:=K$, define
$B(\alpha,\delta)=\{x\in\mathbb{R}^{3};|x-\alpha|<\delta\}$ ($\delta>0$).
Note $K\subset\cup_{\alpha\in K}B(\alpha,\delta)$ and $K$ is compact,
so there exists $\alpha_{1},\cdots,\alpha_{m}\in K$ such that $K\subset\cup_{j=1}^{m}B(\alpha_{j},\delta)$.
Thus, we define 
\[
D_{j,\delta}:=D\cap B(\alpha_{j},\delta)\mbox{ and }D_{\delta}:=\cup_{j=1}^{m}D_{j,\delta}.
\]
It is easy to see that 
\[
\begin{cases}
\int_{D\backslash D_{\delta}}e^{-p\tau(x\cdot\omega-t_{0})A_{t_{0}}^{A}(x')}bdx=O(e^{-pa\tau})\\
\int_{D\backslash D_{\delta}}e^{-p\tau(x\cdot\omega-t_{0})A_{t_{0}}^{B}(x')}bdx=O(e^{-pa\tau})
\end{cases}
\]
where $A_{t_{0}}^{A}(x'),A_{t_{0}}^{B}(x')$ are smooth matrix-valued
functions with bounded entries and their real part strictly greater
than 0. so $\exists a>0$ such that $\mbox{Re}A_{t_{0}}^{A}(x')\geq a>0$
and $\mbox{Re}A_{t_{0}}^{B}(x')\geq a>0$. Let $\alpha_{j}\in K$,
by rotation and translation, we may assume $\alpha_{j}=0$ and the
vector $\alpha_{j}-x_{0}=-x_{0}$ is parallel to $e_{3}=(0,0,1)$.
Therefore, we consider the change of coordinates near each $\alpha_{j}$
as follows:
\[
\begin{cases}
y'=x'\\
y_{3}=x\cdot\rho-t_{0},
\end{cases}
\]
where $x=(x_{1},x_{2},x_{3})=(x',x_{3})$ and $y=(y_{1},y_{2},y_{3})=(y',y_{3})$.
Denote the parametrization of $\partial D$ near $\alpha_{j}$ by
$l_{j}(y')$, then we have the following estimates. Note that the
oscillating-decaying solutions are well-defined in $D$.
\begin{lem}
For $q\leq2$, $\tau\gg1$, we have the following estimates.\\
1. 
\begin{eqnarray*}
\int_{D}|H_{t}(x)|^{q}dx & \le & \tau^{2q-1}\sum_{j=1}^{m}\iint_{|y'|<\delta}e^{-aq\tau l_{j}(y')}dy'+O(\tau^{2q-1}e^{-qa\delta\tau})\\
 &  & +O(\tau^{2q}e^{-qa\tau})+O(\tau e^{-c\tau})+O(\tau^{-2N+5})
\end{eqnarray*}
2. 
\begin{eqnarray*}
\int_{D}|H_{t}|^{2}dx & \geq & C\tau^{3}\sum_{j=1}^{m}\iint_{|y'|<\delta}e^{-2a\tau l_{j}(y')}dy'-C\tau^{3}e^{-2a\delta\tau}\\
 &  & -C\tau e^{-2c\tau}-C\tau^{-2N+5}
\end{eqnarray*}
3. 
\begin{eqnarray*}
\int_{D}|E_{t}(x)|^{q}dx & \le & \tau^{q-1}\sum_{j=1}^{m}\iint_{|y'|<\delta}e^{-aq\tau l_{j}(y')}dy'+O(\tau^{q-1}e^{-qa\delta\tau})\\
 &  & +O(\tau^{q}e^{-qa\tau})+O(\tau^{-1})+O(\tau^{-2N+3})
\end{eqnarray*}
4. 
\begin{eqnarray*}
\int_{D}|E_{t}|^{2}dx & \geq & C\tau\sum_{j=1}^{m}\iint_{|y'|<\delta}e^{-2a\tau l_{j}(y')}dy'-C\tau e^{-2a\delta\tau}\\
 &  & -C\tau^{-1}-C\tau^{-2N+3}
\end{eqnarray*}
\end{lem}
\begin{proof}
The proof is via the representation of the oscillating-decaying solutions
of $(E_{t},H_{t})$. For $\tau\gg1$($\tau\ll\tau^{2}$), we have
\begin{eqnarray*}
\int_{D}|H_{t}|^{q}dx & \leq & C\tau^{2q}\int_{D}e^{-qa\tau(x\cdot\omega-t_{0})}dx+C_{q}\int_{D}|\Gamma_{\chi_{t},b,t,N,\omega}^{A,2}|^{q}dx\\
 &  & +C_{q}\int_{D}|r_{\chi_{t},b,t,N,\omega}^{A,2}|^{q}dx\\
 & \leq & C\tau^{2q}\int_{D_{\delta}}e^{-qa\tau(x\cdot\omega-t_{0})}dx+C\tau^{2q}\int_{D\backslash D_{\delta}}e^{-qa\tau(x\cdot\omega-t_{0})}dx\\
 &  & +C_{q}\int_{D}|\Gamma_{A,B,\gamma,\mu}^{1}|^{q}dx+C_{q}\int_{D}|r_{A,B,\gamma,\mu}^{1}|^{q}dx\\
 & \leq & C\tau^{2q}\sum_{j=1}^{m}\iint_{|y'|<\delta}dy'\int_{l_{j}(y')}^{\delta}e^{-qa\tau y_{3}}dy_{3}+C\tau^{2q}e^{-qa\tau}\\
 &  & +C\|\Gamma_{\chi_{t},b,t,N,\omega}^{A,2}\|_{L^{2}(D)}^{2}+C\|r_{\chi_{t},b,t,N,\omega}^{A,2}\|_{L^{2}(D)}^{2}\\
 & \leq & C\tau^{2q-1}\sum_{j=1}^{m}\iint_{|y'|<\delta}e^{-aq\tau l_{j}(y')}dy'-\dfrac{C}{q}\tau^{2q-1}e^{-qa\delta\tau}\\
 &  & +C\tau^{2q}e^{-qa\tau}+C\tau e^{-ca\tau}+C\tau^{-2N+5},
\end{eqnarray*}
where $c$ is a positive constant and $a$ depending only on $a_{A},a_{B}$.
For the lower bound of $\int_{D}|H_{t}|^{2}dx$, we have 
\begin{eqnarray*}
\int_{D}|H_{t}|^{2}dx & \geq & C\tau^{4}\int_{D}e^{-2a\tau(x\cdot\omega-t_{0})}dx-C\|\Gamma_{\chi_{t},b,t,N,\omega}^{A,2}\|_{L^{2}(\Omega_{t_{0}}(\omega))}^{2}\\
 &  & -C\|r_{\chi_{t},b,t,N,\omega}^{A,2}\|_{L^{2}(\Omega_{t_{0}}(\omega))}^{2}\\
 & \geq & C\tau^{4}\int_{D_{\delta}}e^{-2a\tau(x\cdot\omega-t_{0})}dx-C\tau e^{-c\tau}-C\tau^{-2N+5}.\\
 & \geq & C\tau^{3}\sum_{j=1}^{m}\iint_{|y'|<\delta}e^{-2a\tau l_{j}(y')}dy'-C\tau^{3}e^{-2a\delta\tau}\\
 &  & -C\tau e^{-ca\tau}-C\tau^{-2N+5}.
\end{eqnarray*}
It is similar to prove the remaining case, so we omit the proof.\end{proof}
\begin{lem}
We have the following estimate 
\[
\dfrac{\|H_{t}\|_{L^{2}(D)}^{2}}{\|E_{t}\|_{L^{2}(D)}^{2}}\geq O(\tau^{2}),\mbox{ }\tau\gg1.
\]
\end{lem}
\begin{proof}
Since $\partial D$ is Lipschitz, we have $l_{j}(y')\le C|y'|$. Therefore
we have the following estimate 
\begin{eqnarray*}
C\tau^{3}\sum_{j=1}^{m}\iint_{|y'|<\delta}e^{-2a\tau l_{j}(y')}dy' & \geq & C\tau^{3}\sum_{j=1}^{m}\iint_{|y'|<\delta}e^{-2a\tau|y'|}\\
 & \geq & C\tau\sum_{j=1}^{m}\iint_{|y'|<\tau\delta}e^{-2a|y'|}dy'\\
 & = & O(\tau).
\end{eqnarray*}
Then we use lemma 4.4 to get 
\begin{eqnarray*}
\dfrac{\|H_{t}\|_{L^{2}(D)}^{2}}{\|E_{t}\|_{L^{2}(D)}^{2}} & \geq & C\tau^{2}\dfrac{1-\tfrac{Ce^{-2a\delta\tau}+C\tau^{-2}e^{-2c\tau}+C\tau^{-2N+2}}{\sum_{j=1}^{m}\iint_{|y'|<\delta}e^{-2a\tau l_{j}(y')}dy'}}{1-\tfrac{O(e^{-2\delta a\tau})+O(\tau e^{-ca\tau})+O(\tau^{-2N+2})}{\sum_{j=1}^{m}\iint_{|y'|<\delta}e^{-2a\tau l_{j}(y')}dy'}}\\
 & = & O(\tau^{2})\mbox{ (if }\tau\gg1).
\end{eqnarray*}
\end{proof}
\begin{lem}
If $t=h_{D}(\rho)$, then for some positive constant $C$, we have
\[
\liminf_{\tau\to\infty}\int_{D}\tau|\nabla\times H_{t}|^{2}dx\geq C.
\]
\end{lem}
\begin{proof}
Since $l_{j}(y')\leq C|y'|$, we have 
\begin{eqnarray*}
\int_{D}|\nabla\times H_{t}(x)|^{2}dx & \geq & C\int_{D}|E_{t}(x)|^{2}dx\\
 & \geq & C\tau\sum_{j=1}^{m}\iint_{|y'|<\delta}e^{-2a\tau l_{j}(y')}dy'-C\tau e^{-2a\delta\tau}\\
 &  & -C\tau^{-1}-C\tau^{-2N+3}\\
 & \geq & C\tau\sum_{j=1}^{m}\iint_{|y'|<\delta}e^{-2a\tau|y'|}dy'-C\tau e^{-2a\delta\tau}\\
 &  & -C\tau^{-1}-C\tau^{-2N+3}\\
 & \geq & C\tau[\tau^{-2}\sum_{j=1}^{m}\iint_{|y'|<\tau\delta}e^{-2a|y'|}dy']-C\tau e^{-2a\delta\tau}\\
 &  & -C\tau^{-1}-C\tau^{-2N+3}\mbox{ (as }\tau\gg1).
\end{eqnarray*}
Therefore, we have 
\[
\liminf_{\tau\to\infty}\int_{D}\tau|\nabla\times H_{t}|^{2}dx\geq C.
\]
\end{proof}
\begin{lem}
For $p\in(\max\{\tfrac{4}{3},\tfrac{2+\delta}{1+\delta}\},2]$. we
have the following 
\[
\lim_{\eta\to0}\limsup_{\ell\to\infty}\dfrac{\|\widetilde{H_{\eta,\ell}}\|_{L^{2}(\Omega)}^{2}}{\|\nabla\times H_{t}\|_{L^{2}(D)}^{2}}\leq C\tau^{1-\frac{2}{p}}\mbox{ }(\tau\gg1).
\]
\end{lem}
\begin{proof}
From the proposition 5.2, we have 
\[
\lim_{\eta\to0}\limsup_{\ell\to\infty}\|\widetilde{H_{\eta,\ell}}\|_{L^{2}(\Omega)}\leq C\|\nabla\times H_{t}\|_{L^{p}(D)}.
\]
Then it is easy to see the conclusion.\end{proof}
\begin{rem}
Recall that the sequence $\{H_{\eta,\ell}\}$ converges to $H_{t+\eta}$
in $H(curl,K)$ as $\ell\to\infty$ for all compact subset $D\Subset K\Subset\Omega$
and $H_{t+\eta}\to H_{t}$ in $H^{2}(\Omega_{t}(\omega))$ as $\eta\to0$,
so we have 
\[
\|\nabla\times H_{\eta,\ell}\|_{L^{p}(D)}\to\|\nabla\times H_{t}\|_{L^{p}(D)}\mbox{ and }\|H_{\eta,\ell}\|_{L^{2}(D)}\to\|H_{t}\|_{L^{2}(D)}
\]
as $\ell\to\infty$, $\eta\to0$.
\end{rem}

\subsubsection{End of the proof of Theorem 1.1 for the penetrable case}

First, we prove the case $t<h_{D}(\rho)$. From (\ref{iden1}), we
have
\begin{eqnarray}
-\tau^{-1}I_{\rho}^{\eta,\ell}(\tau,t) & = & \int_{\Omega}\left((\epsilon^{-1}-\epsilon_{0}^{-1})\nabla\times H_{\eta,\ell}\right)\cdot(\nabla\times\overline{H_{\eta,\ell}})dx\nonumber \\
 &  & -\int_{\Omega}(\epsilon^{-1}\nabla\times\widetilde{H_{\eta,\ell}})\cdot(\nabla\times\overline{\widetilde{H_{\eta,\ell}}})dx-k^{2}\int_{\Omega}\mu|\widetilde{H_{\eta,\ell}}|^{2}dx.\label{eq:iden1-1}
\end{eqnarray}
Note that $(\widetilde{E_{\epsilon,\ell}},\widetilde{H_{\epsilon.\ell}})$
satisfies
\[
\begin{cases}
\nabla\times\widetilde{E_{\eta,\ell}}-ik\mu\widetilde{H_{\eta,\ell}}=0 & \mbox{ in }\Omega,\\
\nabla\times\widetilde{H_{\eta,\ell}}+ik\gamma\widetilde{E_{\eta,\ell}}=ik(\epsilon_{0}-\epsilon)E_{\eta,\ell} & \mbox{ in }\Omega,
\end{cases}
\]
and rewrite it as 
\begin{equation}
\nabla\times(\epsilon^{-1}\nabla\times\widetilde{E_{\eta,\ell}})-k^{2}\gamma\widetilde{E_{\eta,\ell}}=k^{2}(\epsilon-\epsilon_{0})E_{\eta,\ell}.\label{eq:iden 1-2}
\end{equation}
Thus, we can use the same argument from the Remark 5.4 again to (\ref{eq:iden 1-2}),
it is easy to see 
\[
\|\widetilde{E_{\eta,\ell}}\|_{L^{2}(\Omega)}\leq C\|E_{\eta,\ell}\|_{L^{2}(D)}.
\]
In addition, we use the Maxwell's equation and $\epsilon-\epsilon_{0}=-\epsilon_{D}\chi_{D}$,
then we have 
\begin{eqnarray}
\int_{\Omega}(\epsilon^{-1}\nabla\times\widetilde{H_{\eta,\ell}})\cdot(\nabla\times\overline{\widetilde{H_{\eta,\ell}}})dx & = & \int_{\Omega}(-ik\epsilon\widetilde{E_{\eta,\ell}}+ik(\epsilon_{0}-\epsilon)E_{\eta,\ell}))\cdot(\nabla\times\overline{\widetilde{H_{\eta,\ell}}})dx\nonumber \\
 & \le & C\int_{\Omega}|\widetilde{E_{\eta,\ell}}|^{2}dx+C\int_{D}|E_{\eta,\ell}|^{2}dx\label{eq:end1}\\
 & \leq & C\int_{D}|E_{\eta,\ell}|^{2}dx.\nonumber 
\end{eqnarray}
Thus, from (\ref{eq:iden1-1}), Proposition 5.3, Lemma 5.5 and (\ref{eq:end1}),
we can obtain
\begin{eqnarray*}
|\dfrac{1}{\tau}I_{\rho}^{\eta,\ell}(\tau,t)| & \leq & \|E_{\eta,\ell}\|_{H(curl,D)}^{2}+\|H_{\eta,\ell}\|_{H(curl,D)}^{2}.
\end{eqnarray*}
From taking $\ell\to\infty$ and $\eta\to0$, we have 
\begin{eqnarray*}
|\dfrac{1}{\tau}I_{\rho}(\tau,t)| & \leq & |\tau\sum_{j=1}^{m}\iint_{|y'|<\delta}e^{-2a\tau l_{j}(y')}dy'+O(\tau^{2}e^{-2a\delta\tau})\\
 &  & +O(\tau^{2}e^{-2a\tau})+O(\tau^{-3})+O(\tau^{-2N+3})\\
 & \leq & O(\tau^{-1})+O(\tau^{2}e^{-2a\delta\tau})\\
 &  & +O(\tau^{2}e^{-2a\tau})+O(\tau^{-3})+O(\tau^{-2N+3}).
\end{eqnarray*}
In particular, we get 
\[
\limsup_{\tau\to\infty}|\dfrac{1}{\tau}I_{\rho}(\tau,t)|=0.
\]

Second, we prove the case $t=h_{D}(\rho)$.\textbf{}\\
\textbf{Case 1}. $\xi\cdot(\gamma^{-1}-\gamma_{0}^{-1})\xi\geq\Lambda|\xi|^{2}$
for all $\xi\in\mathbb{R}^{3}$ for some $\Lambda>0$.

From the inequality in Lemma 5.1, we have 
\begin{eqnarray*}
-\tau^{-1}I_{\rho}^{\eta,\ell} & \geq & \int_{D}[\epsilon(\epsilon-\epsilon_{0}^{-1})^{-1}\epsilon_{0}^{-1}\nabla\times H_{\eta,\ell}]\cdot(\nabla\times\overline{H_{\eta,\ell}})dx-k^{2}\int_{\Omega}\mu|\widetilde{H_{\eta,\ell}}|^{2}dx\\
 &  & -k^{2}\int_{\Omega}\mu|\widetilde{H_{\eta,\ell}}|^{2}dx\\
 & \geq & C\int_{D}|\nabla\times H_{\eta,\ell}|^{2}dx-c\|\widetilde{H_{\eta,\ell}}\|_{L^{2}(\Omega)}^{2}.
\end{eqnarray*}
By using the definition $I_{\rho}(\tau,t):=\lim_{\eta\to0}\lim_{\ell\to\infty}I_{\rho}^{\epsilon,\ell}(\tau,t)$,
$\{H_{\eta,\ell}\}$ converges to $H_{t}$ in $H(curl,K)$ for all
compact subset $D\Subset K\Subset\Omega$ as $\ell\to\infty$, $\eta\to0$,
we have 
\begin{eqnarray*}
\dfrac{-I_{\rho}(\tau,t)}{\|\nabla\times H_{t}\|_{L^{2}(D)}^{2}} & \geq & C\tau\left[1-C\lim_{\epsilon\to0}\limsup_{\ell\to\infty}\dfrac{\|\widetilde{H_{\ell}}\|_{L^{2}(\Omega)}^{2}}{\|\nabla\times H_{t}\|_{L^{2}(D)}^{2}}\right]\\
 & \geq & C\tau(1-C\tau^{1-\frac{2}{p}}).
\end{eqnarray*}
Hence, using Lemma 4.7 we deduce that for $\tau\gg1$, 
\[
|I_{\rho}(\tau,h_{D}(\rho))|\geq C>0
\]
which finishes the proof.\\
\textbf{Case }2. $\xi\cdot(\gamma_{0}^{-1}-\gamma^{-1})\xi\geq\lambda|\xi|^{2}$
for all $\xi\in\mathbb{R}^{3}$ for some $\lambda>0$.

Similarly, using the inequality in Lemma 4.1, we have 
\[
\tau^{-1}I_{\rho}^{\eta,\ell}(\tau,t)\geq\int_{D}((\epsilon_{0}^{-1}-\epsilon^{-1})\nabla\times H_{\eta,\ell})\cdot(\nabla\times\overline{H_{\eta,\ell}})dx-k^{2}\int_{\Omega}\mu|H_{\eta,\ell}|^{2}dx.
\]
Then use the same argument as in \textbf{Case 1} we can finish the
proof.

\subsection{Impenetrable Case}

We give the proof of the second part of Theorem 1.1, since it is the
hardest part. The other cases are easy since we have proved it in
the penetrable case. In addition, the upper bound is easy because
of the well-posedness and the $L^{p}$ estimate for the indicator
function, but the lower bound is not easy to see. In the following
proof, we will use the layer potential properties for the exterior
isotropic Maxwell's equation (with the Silver-M$\ddot{\mathrm{u}}$ller
radiation condition) and the perturbation argument from the anisotropic
Maxwell's equation compared with the isotropic case. In the impenetrable
case, we have chosen the oscillating-decaying solution as the following
form 
\[
\begin{cases}
E_{t}=G_{B}^{2}(x)e^{i\tau x\cdot\xi}e^{-\tau(x\cdot\omega-t)A_{t}^{B}(x')}b+\Gamma_{\chi_{t},b,t,N,\omega}^{B,2}(x,\tau)+r_{\chi_{t},b,t,N,\omega}^{,B,2}(x,\tau),\\
H_{t}=G_{B}^{1}(x)e^{i\tau x\cdot\xi}e^{-\tau(x\cdot\omega-t)A_{t}^{B}(x')}b+\Gamma_{\chi_{t},b,t,N,\omega}^{B,1}(x,\tau)+r_{\chi_{t},b,t,N,\omega}^{B,1}(x,\tau),
\end{cases}
\]
where $G_{B}^{1}(x)=O(\tau)$ and $G_{B}^{2}(x)=O(\tau^{2})$ and
$\Gamma_{\chi_{t},b,t,N,\omega}^{B,j}$ satisfies (\ref{eq:1.7})
for $|\alpha|=j$ and $r_{\chi_{t},b,t,N,\omega}^{B,j}$ satisfies
(\ref{eq:1.7}) for $k=j$.

We start by the following lemma.
\begin{lem}
Assume that $\mu$ is a smooth scalar function and $\gamma$ is a
matrix-valued function. Let $(E,H)\in H(curl;\Omega\backslash\bar{D})\times H(curl;\Omega\backslash\bar{D})$
be a solution of the problem 
\begin{equation}
\begin{cases}
\nabla\times E-ik\mu H=0 & \mbox{ in }\Omega\backslash\bar{D},\\
\nabla\times H+i\epsilon E=0 & \mbox{ in }\Omega\backslash\bar{D},\\
\nu\times E=f & \mbox{ on }\partial\Omega,\\
\nu\times H=0 & \mbox{ on }\partial D,
\end{cases}\label{eq:6.1}
\end{equation}
with $f\in TH^{-1/2}(\partial\Omega)$. If we put $f_{\eta,\ell}=\nu\times E_{\eta,\ell}$
with $\{E_{\eta,\ell}\}$ is obtained by the Runge approximation property.
Then we have the identity
\begin{eqnarray*}
-\dfrac{1}{\tau}I_{\rho}^{\eta,\ell}(\tau,t) & = & -\int_{D}\{|\nabla\times E_{\eta,\ell}(x)|^{2}-k^{2}|E_{\eta,\ell}(x)|^{2}\}dx\\
 &  & -\int_{\Omega\backslash\bar{D}}\{|\nabla\times\widetilde{E_{\eta,\ell}}(x)|^{2}-k^{2}|\widetilde{E_{\eta,\ell}}(x)|^{2}\}dx\\
 & = & \int_{D}\{|\nabla\times H_{\eta,\ell}(x)|^{2}-k^{2}|H_{\eta,\ell}(x)|^{2}\}dx\\
 &  & +\int_{\Omega\backslash\bar{D}}\{|\nabla\times\widetilde{H_{\eta,\ell}}(x)|^{2}-k^{2}|\widetilde{H_{\eta,\ell}}(x)|^{2}\}dx
\end{eqnarray*}
and the inequality 
\[
-\dfrac{1}{\tau}I_{\rho}^{\eta,\ell}(\tau,t)\geq\int_{D}\{|\nabla\times H_{\eta,\ell}(x)|^{2}-k^{2}|H_{\eta,\ell}(x)|^{2}\}dx-k^{2}\int_{\Omega\backslash\bar{D}}|\widetilde{H_{\eta,\ell}}(x)|^{2}\}dx,
\]
where $\widetilde{E_{\eta,\ell}}=E-E_{\eta,\ell}$ and $\widetilde{H_{\eta,\ell}}=H-H_{\eta,\ell}$
are described in section 5.\end{lem}
\begin{proof}
Use the integration by parts and the boundary condition, we have 
\[
\int_{\Omega\backslash\bar{D}}\epsilon^{-1}(\nabla\times E)\cdot\overline{(\nabla\times\widetilde{E_{\eta,\ell}}})-k^{2}\epsilon E\cdot\overline{\widetilde{E_{\eta,\ell}}}dx=-(\int_{\partial\Omega}-\int_{\partial D})ik(\nu\times H)\cdot\overline{\widetilde{E_{\eta,\ell}}}dS=0.
\]
Adding this to 
\begin{eqnarray*}
I_{\rho}^{\eta,\ell} & = & \int_{\partial\Omega}(\nu\times E_{\eta,\ell})\cdot(\overline{-ikH+ikH_{\eta,\ell}})dS\\
 & = & \int_{\Omega\backslash\bar{D}}-(\mu^{-1}\nabla\times E_{\eta,\ell})\cdot(\overline{\nabla\times E})+k^{2}(\mu E_{\eta,\ell})\cdot\bar{E}dx\\
 &  & +\int_{\Omega}\mu^{-1}|\nabla\times E_{\eta,\ell}|^{2}-k^{2}(\mu E_{\eta,\ell})\cdot\overline{E_{\eta,\ell}}dx+\int_{\partial D}(\nu\times E_{\eta,\ell})\cdot(\overline{-ikH})dS
\end{eqnarray*}
due to the zero boundary condition on $\partial D$ we have the last
term is vanishing.
\end{proof}
From the above estimate, it only need to control the lower order term
$\int_{\Omega\backslash\bar{D}}|\widetilde{H_{\eta,\ell}}(x)|^{2}dx$.

\subsubsection{Estimate of the lower order term $\widetilde{H_{\eta,\ell}}$}
\begin{prop}
Let $\Omega$ be a $C^{1}$ domain, $D\Subset\Omega$ be Lipschitz.
Then there exists a positive constant $C$ independent of $(\widetilde{E_{\eta,\ell}},\widetilde{H_{\eta,\ell}})$
and $(E_{\eta,\ell},H_{\eta,\ell})$ such that 
\[
\int_{\Omega\backslash\bar{D}}|\widetilde{H_{\eta,\ell}}(x)|^{2}dx\leq C\{\|\nabla\times H_{\eta,\ell}\|_{L^{p}(D)}^{2}+\|H_{\eta,\ell}\|_{H^{s+1/2}(D)}^{2}\},
\]
for all $p$ and $s$ such that $\max\{2-\delta,4/3\}<p\leq2$ and
$0<s\leq1$ with $\delta>0$.\end{prop}
\begin{proof}
\textbf{Step 1.} Before proving the Proposition 6.2, we consider the
anisotropic Maxwell's equation in $\Omega$ as follows:
\begin{equation}
\begin{cases}
\nabla\times E_{\eta,\ell}-ik\mu H_{\eta,\ell}=0 & \mbox{ in }\Omega,\\
\nabla\times H_{\eta,\ell}+ik\epsilon E_{\eta,\ell}=0 & \mbox{ in }\Omega,\\
\nu\times E_{\eta,\ell}:=f_{\eta,\ell}\in TH^{-1/2}(\partial\Omega) & \mbox{ on }\partial\Omega,
\end{cases}\label{eq:6.2}
\end{equation}
where $E_{\eta,\ell}$ and $H_{\eta,\ell}$ are solutions of the anisotropic
Maxwell's equation. Since $\widetilde{E_{\eta,\ell}}=E-E_{\eta,\ell}$,
$\widetilde{H_{\eta,\ell}}=H-H_{\eta,\ell}$, we have 
\begin{equation}
\begin{cases}
\nabla\times\widetilde{E_{\eta,\ell}}-ik\mu\widetilde{H_{\eta,\ell}}=0 & \mbox{ in }\Omega\backslash\bar{D},\\
\nabla\times\widetilde{H_{\eta,\ell}}+ik\gamma\widetilde{E_{\eta,\ell}}=0 & \mbox{ in }\Omega\backslash\bar{D},\\
\nu\times\widetilde{E_{\eta,\ell}}=0 & \mbox{ on }\partial\Omega,\\
\nu\times\widetilde{H_{\eta,\ell}}=-\nu\times H_{\eta,\ell} & \mbox{ on }\partial D.
\end{cases}\label{eq:6.3}
\end{equation}
\textbf{Step 2.} Let $(E_{\eta,\ell}^{ex},H_{\eta,\ell}^{ex})$ be
the solution of the following well posed exterior Maxwell's problem
\begin{equation}
\begin{cases}
\nabla\times E_{\eta,\ell}^{ex}-ikH_{\eta,\ell}^{ex}=0 & \mbox{ in }\mathbb{R}^{3}\backslash\bar{D},\\
\nabla\times H_{\eta,\ell}^{ex}+ikE_{\eta,\ell}^{ex}=0 & \mbox{ in }\mathbb{R}^{3}\backslash\bar{D},\\
\nu\times H_{\eta,\ell}^{ex}=-\nu\times H_{\eta,\ell} & \mbox{ on }\partial D,\\
E_{\eta,\ell}^{ex},H_{\eta,\ell}^{ex}\mbox{ satisfiy the Silver-M\ensuremath{\ddot{\mathrm{u}}}ller radiation condition}.
\end{cases}\label{eq:6.4}
\end{equation}
We can represent these solutions $E_{\eta,\ell}^{ex}$ and $H_{\eta,\ell}^{ex}$
by the following layer potentials
\begin{eqnarray*}
H_{\eta,\ell}^{ex}(x) & := & \nabla\times\int_{\partial D}\Phi_{k}(x,y)f(y)ds(y),\\
E_{\eta,\ell}^{ex}(x) & := & -\dfrac{1}{ik}\nabla\times H_{\eta,\ell}^{ex}(x),\mbox{ }x\in\mathbb{R}^{3}\backslash\partial D,
\end{eqnarray*}
where $\Phi_{k}(x,y)=-\dfrac{e^{ik|x-y|}}{4\pi|x-y|}$, $x,y\in\mathbb{R}^{3}$,
$x\neq y$, is the fundamental solution of the Helmholtz equation
and $f$ is the density. Now, we follow the arguments in section 2.1
of \cite{KS2014} and use the same argument for the isotropic Maxwell's
equation (\ref{eq:6.4}), then we have 
\begin{equation}
\begin{cases}
\|E_{\eta,\ell}^{ex}\|_{L^{p}(\Omega\backslash\bar{D})}\leq C\{\|\nu\times H_{\eta,\ell}\|_{L^{p}(\partial D)}+\|\nabla\times H_{\eta,\ell}\|_{L^{p}(D)}\},\\
\|H_{\eta,\ell}^{ex}\|_{L^{2}(\Omega\backslash\bar{D})}\leq C\{\|\nu\times H_{\eta,\ell}\|_{L^{p}(\partial D)}+\|\nabla\times H_{\eta,\ell}\|_{L^{p}(D)}\},
\end{cases}\label{eq:6.5}
\end{equation}
for $p\in(\dfrac{4}{3},2]$. Moreover, if we define $\mathcal{E}_{\eta,\ell}=\widetilde{E_{\eta,\ell}}-E_{\eta,\ell}^{ex}$,
$\mathcal{H}_{\eta,\ell}=\widetilde{H_{\eta,\ell}}-H_{\eta,\ell}^{ex}$,
then $\mathcal{E}_{\eta,\ell}$ and $\mathcal{H}_{\eta,\ell}$ satisfy
the following Maxwell's equation
\begin{equation}
\begin{cases}
\nabla\times\mathcal{E}_{\eta,\ell}-ik\mu\mathcal{H}_{\eta,\ell}=ik(1-\mu)H_{\eta,\ell}^{ex} & \mbox{ in }\Omega\backslash\bar{D},\\
\nabla\times\mathcal{H}_{\eta,\ell}+ik\epsilon\mathcal{E}_{\eta,\ell}=ik(\gamma-I_{3})E_{\eta,\ell}^{ex} & \mbox{ in }\Omega\backslash\bar{D},\\
\nu\times\mathcal{H}_{\eta,\ell}=0 & \mbox{ on }\partial\Omega,\\
\nu\times\mathcal{E}_{\eta,\ell}=-\nu\times E_{\eta,\ell}^{ex} & \mbox{ on }\partial D.
\end{cases}\label{eq:6.6}
\end{equation}
\textbf{Step 3.} Now we decompose $\mathcal{E}_{\eta,\ell}=\mathcal{E}_{\eta,\ell}^{1}+\mathcal{E}_{\eta,\ell}^{2}$
and $\mathcal{H}_{\eta,\ell}=\mathcal{H}_{\eta,\ell}^{1}+\mathcal{H}_{\eta,\ell}^{2}$,
where $(\mathcal{E}_{\eta,\ell}^{1},\mathcal{H}_{\eta,\ell}^{1})$
satisfies the following zero boundary Maxwell's equation 
\begin{equation}
\begin{cases}
\nabla\times\mathcal{E}_{\eta,\ell}^{1}-ik\mu\mathcal{H}_{\eta,\ell}^{1}=ik(1-\mu)H_{\eta,\ell}^{ex} & \mbox{ in }\Omega\backslash\bar{D},\\
\nabla\times\mathcal{H}_{\eta,\ell}^{1}+ik\epsilon\mathcal{E}_{\eta,\ell}^{1}=ik(\epsilon-I_{3})E_{\eta,\ell}^{ex} & \mbox{ in \ensuremath{\Omega}}\backslash\bar{D},\\
\nu\times\mathcal{E}_{\eta,\ell}^{1}=\nu\times\mathcal{H}_{\eta,\ell}^{1}=0 & \mbox{ on }\partial(\Omega\backslash\bar{D}),
\end{cases}\label{eq:6.7}
\end{equation}
and $(\mathcal{E}_{\eta,\ell}^{2},\mathcal{H}_{\eta,\ell}^{2})$ satisfies
\begin{equation}
\begin{cases}
\nabla\times\mathcal{E}_{\eta,\ell}^{2}-ik\mu\mathcal{H}_{\eta,\ell}^{2}=0 & \mbox{ in }\Omega\backslash\bar{D},\\
\nabla\times\mathcal{H}_{\eta,\ell}^{2}+ik\gamma\mathcal{E}_{\eta,\ell}^{2}=0 & \mbox{ in }\Omega\backslash\bar{D},\\
\nu\times\mathcal{H}_{\eta,\ell}^{2}=0 & \mbox{ on }\partial\Omega,\\
\nu\times\mathcal{E}_{\eta,\ell}^{2}=-\nu\times E_{\eta,\ell}^{ex} & \mbox{ on }\partial D.
\end{cases}\label{eq:6.8}
\end{equation}

First, we deal with the equation (\ref{eq:6.7}) by using the $L^{p}$
estimate in $\Omega\backslash\bar{D}$. Note that $(\mathcal{E}_{\eta,\ell}^{1},\mathcal{H}_{\eta,\ell}^{1})$
satisfies (\ref{eq:6.7}), then we have 
\[
\begin{cases}
\nabla\times(\epsilon^{-1}\nabla\times\mathcal{E}_{\eta,\ell}^{1})-k^{2}\gamma\mathcal{E}_{\eta,\ell}^{1}=ik\nabla\times[(\mu^{-1}-1)H_{\eta,\ell}^{ex}]+ik(\gamma-I_{3})E_{\eta,\ell}^{ex} & \mbox{ in }\Omega\backslash\bar{D},\\
\nu\times\mathcal{E}_{\eta,\ell}^{1}=0 & \mbox{ on }\partial(\Omega\backslash\bar{D}),
\end{cases}
\]
and 
\[
\begin{cases}
\nabla\times(\epsilon^{-1}\nabla\times\mathcal{H}_{\eta,\ell}^{1})-k^{2}\mu\mathcal{H}_{\eta,\ell}^{1}=ik\nabla\times[(I_{3}-\epsilon^{-1})E_{\eta,\ell}^{ex}]+ik(1-\mu)H_{\eta,\ell}^{ex} & \mbox{ in }\Omega\backslash\bar{D},\\
\nu\times\mathcal{H}_{\eta,\ell}^{1}=0 & \mbox{ on }\partial(\Omega\backslash\bar{D}).
\end{cases}
\]
Now, if we use the same method in the proof of the Proposition 5.3,
we will obtain 
\begin{equation}
\begin{cases}
\|\mathcal{E}_{\eta,\ell}^{1}\|_{L^{p}(\Omega\backslash\bar{D})}+\|\nabla\times\mathcal{E}_{\eta,\ell}^{1}\|_{L^{p}(\Omega\backslash\bar{D})}\leq C\{\|H_{\eta,\ell}^{ex}\|_{L^{p}(\Omega\backslash\bar{D})}+\|E_{\epsilon,\ell}^{ex}\|_{L^{2}(\Omega\backslash\bar{D})}\},\\
\|\mathcal{H}_{\eta,\ell}^{1}\|_{L^{p}(\Omega\backslash\bar{D})}+\|\nabla\times\mathcal{H}_{\eta,\ell}^{1}\|_{L^{P}(\Omega\backslash\bar{D})}\leq C\{\|E_{\eta,\ell}^{ex}\|_{L^{p}(\Omega\backslash\bar{D})}+\|H_{\eta,\ell}^{ex}\|_{L^{2}(\Omega\backslash\bar{D})}.
\end{cases}\label{eq:6.9}
\end{equation}
for any $\dfrac{4}{3}<p\leq2$. If we combine (\ref{eq:6.5}) and
(\ref{eq:6.9}) together, we have 
\begin{equation}
\|H_{\eta,\ell}^{1}\|_{L^{p}(\Omega\backslash\bar{D})}\leq C\{\|\nu\times H_{\eta,\ell}\|_{L^{p}(\partial D)}+\|\nabla\times H_{\eta,\ell}\|_{L^{p}(D)}\}.\label{eq:6.10}
\end{equation}

For $(\mathcal{E}_{\eta,\ell}^{2},\mathcal{H}_{\eta,\ell}^{2})$,
we apply the $L^{2}$-theory for the anisotropic Maxwell's equation,
we get 
\[
\|\mathcal{H}_{\eta,\ell}^{2}\|_{L^{2}(\Omega\backslash\bar{D})}\leq\|\mathcal{E}_{\eta,\ell}^{2}\|_{H(curl,\Omega\backslash\bar{D})}\leq C\|\nu\times\mathcal{E}_{\eta,\ell}^{2}\|_{H^{-1/2}(\partial\Omega)}\leq C\|\nu\times E_{\eta,\ell}^{ex}\|_{H^{-1/2}(\partial\Omega)}.
\]
Moreover, following the proof in the Lemma 2.3 of \cite{KS2014},
we have 
\[
\|\nu\times E_{\eta,\ell}^{ex}\|_{H^{-1/2}(\partial\Omega)}\leq C\|f\|_{L^{p}(\partial D)},\mbox{ }\forall p\geq1,
\]
and 
\begin{equation}
\|\mathcal{H}_{\eta,\ell}^{2}\|_{L^{2}(\Omega\backslash\bar{D})}\leq C\{\|\nu\times H_{\eta,\ell}\|_{L^{p}(\partial D)}^{2}+\|\nabla\times H_{\eta,\ell}\|_{L^{p}(D)}^{2}\},\label{eq:6.11}
\end{equation}
for all $p\in(\dfrac{4}{3},2]$. Recall that $\mathcal{H}_{\eta,\ell}=\mathcal{H}_{\eta,\ell}^{1}+\mathcal{H}_{\eta,\ell}^{2}$,
by using (\ref{eq:6.10}) and (\ref{eq:6.11}), then we have 
\begin{equation}
\|\mathcal{H}_{\eta,\ell}\|_{L^{2}(\Omega\backslash\bar{D})}\leq C\{\|\nu\times H_{\eta,\ell}\|_{L^{p}(\partial D)}+\|\nabla\times H_{\eta,\ell}\|_{L^{p}(D)}\}\label{eq:6.12}
\end{equation}
for all $p\in(\dfrac{4}{3},2]$. Combining (\ref{eq:6.5}), (\ref{eq:6.12})
and $\widetilde{H_{\eta,\ell}}=\mathcal{H}_{\eta,\ell}+H_{\eta,\ell}^{ex}$,
we get 
\begin{eqnarray}
\int_{\Omega\backslash\bar{D}}|\widetilde{H_{\eta,\ell}}(x)|^{2}dx & \leq & \|\mathcal{H}_{\eta,\ell}\|_{L^{2}(\Omega\backslash\bar{D})}+\|H_{\eta,\ell}^{ex}\|_{L^{2}(\Omega\backslash\bar{D})}\nonumber \\
 & \leq & C\{\|\nu\times H_{\eta,\ell}\|_{L^{p}(\partial D)}^{2}+\|\nabla\times H_{\eta,\ell}\|_{L^{p}(D)}^{2}\}\label{eq:6.13}
\end{eqnarray}
for all $p\in(\dfrac{4}{3},2]$. Finally, for $s>0$ and $p\leq2$
we have $H^{s}(\partial D)\subset L^{2}(\partial D)\subset L^{p}(\partial D)$,
then we reduce that 
\[
\|\nu\times H_{\eta,\ell}\|_{L^{p}(\partial D)}\leq C\|H_{\eta,\ell}\|_{L^{p}(\partial D)}\leq C\|H_{\eta,\ell}\|_{H^{s}(\partial D)}.
\]
Note that the trace map from $H^{s+1/2}(D)\to H^{s}(\partial D)$
is bounded for all $0<s\leq1$. So the estimate (\ref{eq:6.13}) will
become
\[
\int_{\Omega\backslash\bar{D}}|\widetilde{H_{\eta,\ell}}(x)|^{2}dx\leq C\{\|H_{\eta,\ell}\|_{H^{s+1/2}(D)}^{2}+\|\nabla\times H_{\eta,\ell}\|_{L^{p}(D)}^{2}\},
\]
for all $p\in(\dfrac{4}{3},2]$ and $0<s\leq1$.\end{proof}
\begin{rem}
Now, if we take $\ell\to\infty$ and $\epsilon\to0$, we will get
\[
\lim_{\eta\to0}\limsup_{\ell\to\infty}\int_{\Omega\backslash\bar{D}}|\widetilde{H_{\eta,\ell}}(x)|^{2}dx\leq C\{\|H_{t}\|_{H^{s+1/2}(D)}^{2}+\|\nabla\times H_{t}\|_{L^{p}(D)}^{2}\},
\]
where $H_{t}$ is the oscillating-decaying solution defined on $\Omega_{t}(\omega)$.
\end{rem}
We have the following lemmas for the oscillating-decaying solutions
in the same way as we did in section 5, so we omit the proofs.
\begin{lem}
For $1\leq q<\infty$, $\tau\gg1$, we have the following estimates.\\
1. 
\begin{eqnarray*}
\int_{D}|H_{t}(x)|^{q}dx & \le & \tau^{q-1}\sum_{j=1}^{m}\iint_{|y'|<\delta}e^{-aq\tau l_{j}(y')}dy'+O(\tau^{q-1}e^{-qa\delta\tau})\\
 &  & +O(\tau^{q}e^{-qa\tau})+O(\tau^{-1})+O(\tau^{-2N+3})
\end{eqnarray*}
2. 
\begin{eqnarray*}
\int_{D}|H_{t}|^{2}dx & \geq & C\tau\sum_{j=1}^{m}\iint_{|y'|<\delta}e^{-2a\tau l_{j}(y')}dy'-C\tau e^{-2a\delta\tau}\\
 &  & -C\tau^{-1}-C\tau^{-2N+3}
\end{eqnarray*}
3. 
\begin{eqnarray*}
\int_{D}|\nabla\times H_{t}(x)|^{q}dx & \le & \tau^{2q-1}\sum_{j=1}^{m}\iint_{|y'|<\delta}e^{-aq\tau l_{j}(y')}dy'+O(\tau^{2q-1}e^{-qa\delta\tau})\\
 &  & +O(\tau^{2q}e^{-qa\tau})+O(\tau e^{-c\tau})+O(\tau^{-2N+5})
\end{eqnarray*}
4. 
\begin{eqnarray*}
\int_{D}|\nabla\times H_{t}(x)|^{2}dx & \geq & C\tau^{3}\sum_{j=1}^{m}\iint_{|y'|<\delta}e^{-2a\tau l_{j}(y')}dy'-C\tau^{3}e^{-2a\delta\tau}\\
 &  & -C\tau e^{-c\tau}-C\tau^{-2N+5}
\end{eqnarray*}

\end{lem}

\begin{lem}
We have the following estimate 
\[
\dfrac{\|H_{t}\|_{L^{2}(D)}^{2}}{\|\nabla\times H_{t}\|_{L^{2}(D)}^{2}}\leq O(\tau^{-2}),\mbox{ }\tau\gg1.
\]

\end{lem}

For $p<2$, we have the following estimate 
\[
\dfrac{\|\nabla\times H_{t}\|_{L^{p}(D)}^{2}}{\|\nabla\times H_{t}\|_{L^{p}(D)}^{2}}\leq C\tau^{1-\frac{2}{p}},\mbox{ }\tau\gg1.
\]

\begin{lem}
If $t=h_{D}(\rho),$then for some positive constant $C,$ we have
\[
\liminf_{\tau\to\infty}\int_{D}\tau|\nabla\times H_{t}|^{2}dx\geq C.
\]

\end{lem}

\subsubsection{End of the proof of Theorem 1.1 for the impenetrable case}

By using the same argument in the penetrable case, it is easy to see
that 
\[
\limsup_{\tau\to\infty}|\dfrac{1}{\tau}I_{\rho}(\tau,t)|=0
\]
for $t>h_{D}(\rho)$. Recall that from Lemma 6.1, we have 
\begin{equation}
-\dfrac{1}{\tau}I_{\rho}^{\eta,\ell}(\tau,t)\geq\int_{D}\{|\nabla\times H_{\eta,\ell}(x)|^{2}-k^{2}|H_{\eta,\ell}(x)|^{2}\}dx-k^{2}\int_{\Omega\backslash\bar{D}}|\widetilde{H_{\eta,\ell}}(x)|^{2}\}dx.\label{eq:6.0}
\end{equation}
By using Proposition 6.2, we deduce 
\[
-\dfrac{1}{\tau}I_{\rho}^{\eta,\ell}(\tau,t)\geq\int_{D}\{|\nabla\times H_{\eta,\ell}(x)|^{2}-k^{2}|H_{\eta,\ell}(x)|^{2}\}dx-C\{\|H_{t}\|_{H^{s+1/2}(D)}^{2}+\|\nabla\times H_{t}\|_{L^{p}(D)}^{2}\},
\]
where $0<s\leq1$ and $\dfrac{4}{3}<p\leq2$. We want to estimate
$\dfrac{\|H_{t}\|_{H^{s+1/2}(D)}^{2}}{\|\nabla\times H_{t}\|_{L^{2}(D)}^{2}}$,
for $0<s\leq1$. Set $r=s+1/2$, then we need to estimate 
\[
\dfrac{\|H_{t}\|_{H^{r}(D)}^{2}}{\|\nabla\times H_{t}\|_{L^{2}(D)}^{2}}
\]
for $r\in(\dfrac{1}{2},\dfrac{3}{2}]$. Using the interpolation inequality,
we have 
\[
\|H_{t}\|_{H^{r}(D)}\leq C\|H_{t}\|_{L^{2}(D)}^{1-r}\|H_{t}\|_{H^{1}(D)}^{r},\mbox{ }0\leq r\leq1.
\]
By the Young's inequality $ab\leq\delta^{-\alpha}\frac{a^{\alpha}}{\alpha}+\delta^{\beta}\frac{b^{\beta}}{\beta}$,
$\frac{1}{\alpha}+\frac{1}{\beta}=1$, we obtain 
\begin{eqnarray}
\|H_{t}\|_{H^{r}(D)}^{2} & \leq & C\left[\dfrac{\delta^{-\alpha}}{\alpha}\|H_{t}\|_{L^{2}(D)}^{2}+\dfrac{\delta^{\beta}}{\beta}\|H_{t}\|_{H^{1}(D)}^{2}\right]\nonumber \\
 & \leq & C\left[\{(1-r)\delta^{-(1-r)^{-1}}+r\delta^{r^{-1}}\}\|H_{0}\|_{L^{2}(D)}^{2}+r\delta^{r^{-1}}\|\nabla H_{t}\|_{L^{2}(D)}^{2}\right].\label{eq:6.14}
\end{eqnarray}
Recall that $H_{t}=G_{B}^{1}(x)e^{i\tau x\cdot\xi}e^{-\tau(x\cdot\omega-t)A_{t}^{B}(x')}b+\Gamma_{\chi_{t},b,t,N,\omega}^{B,1}(x,\tau)+r_{\chi_{t},b,t,N,\omega}^{B,1}(x,\tau)$
is a smooth function with $G_{B}^{1}(x)=O(\tau)$ and $\Gamma_{\chi_{t},b,t,N,\omega}^{B,1}$
satisfies (\ref{eq:1.7}) for $|\alpha|=1$ and $r_{\chi_{t},b,t,N,\omega}^{B,1}$
satisfies (\ref{eq:1.7}) for $k=1$. If we can differentiate $H_{t}$
componentwisely, we will get $\dfrac{\partial H_{t}}{\partial x_{j}}=\dfrac{\partial G_{B}^{1}e^{i\tau x\cdot\xi}e^{-\tau(x\cdot\omega-t)A_{t}^{B}}b}{\partial x_{j}}+\dfrac{\partial\Gamma_{\chi_{t},b,t,N,\omega}^{B,1}}{\partial x_{j}}+\dfrac{\partial r_{\chi_{t},b,t,N,\omega}^{B,1}}{\partial x_{j}}$
and 
\[
\begin{cases}
\|\dfrac{\partial N_{A,B,\gamma,\mu}^{t}}{\partial x_{j}}\|_{L^{2}(D)}^{2}\leq C\tau^{4}\int_{D}e^{-2a\tau(x\cdot\rho-t)}dx,\\
\|\dfrac{\partial\Gamma_{A,B,\gamma,\mu}^{2,t}}{\partial x_{j}}\|_{L^{2}(D)}\leq c\tau^{-1/2}e^{-c\tau}.\\
\|\dfrac{\partial r_{A,B,\gamma,\mu}^{2,t}}{\partial x_{j}}\|_{L^{2}(D)}\leq c\tau^{-N+3/2}.
\end{cases}
\]
 Then by using the same method as before, it is easy to see that 
\begin{eqnarray*}
\|\nabla H_{t}\|_{L^{2}(D)}^{2} & = & \sum_{j=1}^{3}\|\dfrac{\partial H_{t}}{\partial x_{j}}\|_{L^{2}(D)}^{2}\\
 & \leq & C\tau^{4}\int_{D}e^{-2(x\cdot\rho-t)}dx+c\tau^{-1}e^{-2c\tau}+c\tau^{-2N+3}.
\end{eqnarray*}
For $t=h_{D}(\rho)$, we have 
\begin{eqnarray}
\|\nabla H_{t}\|_{L^{2}(D)}^{2} & \leq & C\tau^{4}\int_{D}e^{-2a(x\cdot\rho-h_{D}(\rho))}dx+c\tau^{-1}e^{-2c\tau}+c\tau^{-2N+3}\nonumber \\
 & \leq & C\tau^{4}(\int_{D_{\delta}}+\int_{D\backslash D_{\delta}})e^{-2a(x\cdot\rho-h_{D}(\rho))}dx+c\tau^{-1}e^{-2\tau(s-t)a}\nonumber \\
 &  & +c\tau^{-2N+3}\nonumber \\
 & \leq & C\tau^{4}\sum_{j=1}^{m}\iint_{|y'|<\delta}dy'\int_{l_{j}(y')}^{\delta}e^{-2a\tau y_{3}}dy_{3}+C\tau^{4}e^{-2ac\tau}\nonumber \\
 &  & +c\tau^{-1}e^{-2c\tau}+c\tau^{-2N+3}\nonumber \\
 & \leq & C\tau^{3}\sum_{j=1}^{m}\iint_{|y'|<\delta}e^{-2a\tau l_{j}(y')}dy'-C\tau^{3}e^{-2a\delta\tau}\nonumber \\
 &  & +C\tau^{3}e^{-2ac\tau}+c\tau^{-1}e^{-2c\tau}+c\tau^{-2N+3}.\label{eq:6.15}
\end{eqnarray}
From Lemma 6.4 and (\ref{eq:6.15}), we have 
\begin{equation}
\dfrac{\|\nabla H_{t}\|_{L^{2}(D)}^{2}}{\|\nabla\times H_{t}\|_{L^{2}(D)}^{2}}\leq C.\label{eq:6.16}
\end{equation}
Combining Lemma 6.4, (\ref{eq:6.14}) and (\ref{eq:6.16}) we obtain
\begin{eqnarray*}
\dfrac{\|H_{t}\|_{H^{r}(D)}^{2}}{\|\nabla\times H_{t}\|_{L^{2}(D)}^{2}} & \leq & C\{(1-r)\delta^{-(1-r)^{-1}}+r\delta^{r^{-1}}\}\dfrac{\|H_{t}\|_{L^{2}(D)}^{2}}{\|\nabla\times H_{t}\|_{L^{2}(D)}^{2}}\\
 &  & +Cr\delta^{r^{-1}}\dfrac{\|\nabla H_{t}\|_{L^{2}(D)}^{2}}{\|\nabla\times H_{t}\|_{L^{2}(D)}^{2}}\\
 & \leq & C\{(1-r)\delta^{-(1-r)^{-1}}+r\delta^{r^{-1}}\}O(\tau^{-2})+Cr\delta^{r^{-1}}.
\end{eqnarray*}
We now choose $p\in(\frac{4}{3},2)$, combining (\ref{eq:6.0}), (\ref{eq:6.14})
and (\ref{eq:6.16}) we have 
\begin{eqnarray*}
\dfrac{-\dfrac{1}{\tau}I_{\rho}^{\epsilon,\ell}(\tau,t)}{\|\nabla\times H_{t}\|_{L^{2}(D)}^{2}} & \geq & C-c_{1}\dfrac{\|H_{t}\|_{L^{2}(D)}^{2}}{\|\nabla\times H_{t}\|_{L^{2}(D)}^{2}}-c_{2}\dfrac{\|H_{t}\|_{H^{r}(D)}^{2}}{\|\nabla\times H_{t}\|_{L^{2}(D)}^{2}}-c_{3}\dfrac{\|\nabla\times H_{t}\|_{L^{p}(D)}^{2}}{\|\nabla\times H_{t}\|_{L^{2}(D)}^{2}}\\
 & \geq & C-c_{1}\{(1-r)\delta^{-(1-r)^{-1}}+r\delta^{r^{-1}}\}O(\tau^{-2})-Cr\delta^{r^{-1}}-c_{3}\tau^{1-\frac{2}{p}}\\
 & \geq & C-c_{2}r\delta^{r^{-1}},\mbox{ }\frac{1}{2}<r<1,\mbox{ }\tau\gg1.
\end{eqnarray*}
Hence from Lemma 6.6, we have 
\[
\liminf_{\tau\to\infty}|I_{\rho}(\tau,h_{D}(\rho))|\geq c>0.
\]

\section{Appendix}

\subsection{Construction of the oscillating-decaying solutions $A$ and $B$}

In this subsection, we show how the scheme in \cite{NUW2005(ODS)}
can be used to derive the oscillating-decaying solutions $A$ and
$B$. Recall that $E$ and $H$ satisfy equation (\ref{eq:1.2}),
therefore we need to derive estimates of the higher derivatives for
$A$ and $B$. Note that the main term of $w_{\chi_{t},b,t,N,\omega}^{A}$
(resp. $w_{\chi_{t},b,t,N,\omega}^{B}$) is $\chi_{t}(x')Q_{t}e^{i\tau x\cdot\xi}e^{-\tau(x\cdot\omega-t)A_{t}^{A}(x')}b$
(resp. $\chi_{t}(x')Q_{t}e^{i\tau x\cdot\xi}e^{-\tau(x\cdot\omega-t)A_{t}^{B}(x')}b$),
which can be directly differentiated term by term since it is a multiplication
of smooth functions. So we can calculate $E$ and $H$ directly. For
convenience, we denote $w=w_{\chi_{t},b,t,N,\omega}$ $\gamma=\gamma_{\chi_{t},b,t,N,\omega}(x,\tau)$.
Without loss of generality, we can use the change of coordinates to
assume $t=0$, $\omega=(0,0,1)$ and $\eta=(1,0,0)$, $\zeta=(0,1,0)$.
Define 
\[
\widetilde{Q_{A}}:=e^{-i\tau x'\cdot\xi'}L_{A}(e^{i\tau x'\cdot\xi'}\cdot),\mbox{ }\widetilde{Q_{B}}:=e^{-i\tau x'\cdot\xi'}L_{B}(e^{i\tau x'\cdot\xi'}\cdot)
\]
where $x'=(x_{1},x_{2})$, $\xi'=(\xi_{1},\xi_{2})$ with $|\xi'|=1$
and $L_{A},L_{B}$ have been defined by (\ref{eq:1.4}) and (\ref{eq:.1.5}).
In the following, we will give all the details for the higher derivatives
of $E$ and $H$.

In \cite{NUW2005(ODS)}, the authors used the phase plane method to
get a first order ODE system and we want to decouple the equation
in order to solve it by direct calculations. The method of construction
the oscillating-decaying solution is decomposed into several steps:\\
 \textbf{Step 1.} As mentioned before, we set $\widetilde{Q_{A}}=e^{-i\tau x'\cdot\xi'}L_{A}(e^{i\tau x'\cdot\xi'}\cdot)$,
$\mbox{ }\widetilde{Q_{B}}:=e^{-i\tau x'\cdot\xi'}L_{B}(e^{i\tau x'\cdot\xi'}\cdot)$
and solve $\widetilde{Q_{A}}v_{A}=0$, $\widetilde{Q_{B}}v_{B}=0$.
In the following calculations, we only need to consider $\widetilde{Q_{A}}v_{A}=0$
since $\widetilde{Q_{B}}v_{B}=0$ will follow the same calculations.
Let $Q_{A}=C_{A}\widetilde{Q_{A}}$ be the operator which satisfies
the leading coefficient of $\partial_{3}^{2}$ is $1$ and the existence
of $C_{A}$ is given by the strong ellipticity of $L_{A}$ and we
need to solve $Q_{A}v_{A}=0$ (the same reason for the operator $\widetilde{Q_{B}}$
and $Q_{B}$). Now, We introduce the concept of the order in the following
manner. We consider $\tau,\partial_{3}$ are of order 1, $\partial_{1},\partial_{2}$
are of order 0 and $x_{3}$ is of order $-1$.\\
 \textbf{Step 2.} Use the Taylor expansion with respect to $x_{3}$,
we have 
\begin{eqnarray*}
Q_{A}(x',x_{3}) & = & Q_{A}(x',0)+\cdots+\dfrac{x_{3}^{N-1}}{(N-1)!}\partial_{3}^{N-1}Q_{A}(x',0)+R\\
 & = & Q_{A}^{2}+Q_{A}^{1}+\cdots+Q_{A}^{-N+1}+R
\end{eqnarray*}
where ord$(Q_{A}^{j})=j$ and ord$(R)=-N$. Since we hope that $Q_{A}v_{A}=0$,
we have 
\[
Q_{A}^{2}v_{A}=-(Q_{A}^{2}+Q_{A}^{1}+\cdots+Q_{A}^{-N+1}+R)v_{A}:=f.
\]
\\
 \textbf{Step 3.} Following the paper \cite{NUW2005(ODS)}, we denote
$D_{3}=-i\partial_{3}$, $\rho=(\xi_{1},\xi_{2},0)$ and $\left\langle a,b\right\rangle =(\left\langle a,b\right\rangle _{ik})$
for $a=(a_{1},a_{2},a_{3})$ and $b=(b_{1},b_{2},b_{3})$, where $\left\langle a,b\right\rangle _{ik}=\sum_{jl}C_{ijkl}^{A}a_{j}b_{l}$
with $C_{ijkl}^{A}$ being the leading coefficient of the second order
strongly elliptic operator $L_{A}$. If we set $W=\left[\begin{array}{c}
w_{1}\\
w_{2}
\end{array}\right]$, where 
\[
\begin{cases}
w_{1}=v_{A}\\
w_{2}=-\tau^{-1}\left\langle e_{3},e_{3}\right\rangle _{x_{3}=0}D_{3}v_{A}-\left\langle e_{3},\rho\right\rangle _{x_{3}=0}v_{A}
\end{cases},
\]
and use $f=-(Q_{A}^{2}+Q_{A}^{1}+\cdots+Q_{A}^{-N+1}+R)v_{A}$, then
$W$ will satisfy 
\begin{eqnarray*}
D_{3}W & = & \tau K^{A}W+\left[\begin{array}{c}
0\\
\tau^{-1}\left\langle e_{3},e_{3}\right\rangle _{x_{3}=0}f
\end{array}\right]\\
 & = & (\tau K^{A}+K_{0}^{A}+\cdots+K_{-N}^{A}+S)W
\end{eqnarray*}
where $K^{A}$ is a matrix in depending of $x_{3}$ which can be diagonlizable
by the property of the strong ellipticity of $L_{A}$. Note that each
$K_{j}^{A}$'s only involves the $x'$ derivatives with ord$(K_{j}^{A})=j$,
ord$(S)=-N-1$. It is worth to mention that with the help of such
special $W$, then we can solve the ODE system explicitly.\\
 \textbf{Step 4.} Decompose $K^{A}$ such that 
\[
\widetilde{K^{A}}=\tilde{Q}^{-1}K^{A}\tilde{Q}=\left[\begin{array}{cc}
\widetilde{K_{+}^{A}} & 0\\
0 & \widetilde{K_{-}^{A}}
\end{array}\right],
\]
where spec$(\widetilde{K_{\pm}^{A}})\subset\mathbb{C}_{\pm}:=\{\pm Im\lambda>0\}$
(the existence of $\widetilde{K^{A}}$ and $\tilde{Q}$ were showed
in \cite{NUW2005(ODS)}). If we set $\widehat{W}=\tilde{Q}^{-1}W$,
then 
\[
D_{3}\widehat{W}=(\tau\widetilde{K^{A}}+\widehat{K}_{0}+\cdots+\widehat{K}_{-N}+\widehat{S})\widehat{W},
\]
\\
 \textbf{Step 5.} If we write $\widehat{W}=(I+x_{3}A^{(0)}+B^{(0)})\widetilde{W}^{(0)}$
with $A^{(0)},B^{(0)}$ being differential operators in $\partial_{x'}$
(their coefficients independent of $x_{3}$), then 
\begin{eqnarray*}
D_{3}\widetilde{W}^{(0)} & = & \{\tau\widetilde{K^{A}}+(\widehat{K}_{0}-\tau x_{3}A^{(0)}\widetilde{K^{A}}+\tau x_{3}\widetilde{K^{A}}A^{(0)}-B^{(0)}\widetilde{K^{A}}\\
 &  & +\widetilde{K^{A}}B^{(0)}+iA^{(0)})+\widehat{K'}_{-1}+\cdots\}\widetilde{W}^{(0)}\\
 & := & (\tau\widetilde{K^{A}}+\widetilde{K}_{0}+\widehat{K'}_{-1}+\cdots)\widetilde{W}^{(0)}
\end{eqnarray*}
where ord$(\widehat{K'}_{-1})=-1$ and the remainders are at most
$-2$. We choose $A^{(0)},B^{(0)}$ to be suitable operators and use
the same calculations in \cite{NUW2005(ODS)}, then we will get 
\[
\widetilde{K}_{0}=\left[\begin{array}{cc}
\widetilde{K}_{0}(1,1) & 0\\
0 & \widetilde{K}_{0}(2,2)
\end{array}\right]
\]
to be a diagonal form (here we omit all the details).\\
 \textbf{Step 6.} Finally, following step 5, we can write 
\begin{eqnarray*}
\widehat{W} & = & (I+x_{3}A^{(0)}+\tau^{-1}B^{(0)})(I+x_{3}^{2}A^{(1)}+\tau^{-1}x_{3}B^{(1)}+\tau^{-2}C^{(1)})\cdots\\
 &  & \times(I+x_{3}^{N+1}A^{(N)}+\tau^{-1}x_{3}^{N}B^{(N)}+\tau^{-2}x_{3}^{N-1}C^{(N)})\widetilde{W}^{(N)}
\end{eqnarray*}
with suitable $A^{(j)},B^{(j)}$ and $C^{(j)}$ for $j=0,1,2,\cdots,N$
($C^{(0)}=0$), then $\widetilde{W}^{(N)}$ satisfies 
\[
D_{3}\widetilde{W}^{(N)}=\{\tau\widetilde{K^{A}}+\widetilde{K}_{0}+\cdots+\widetilde{K}_{-N}+\widetilde{S}\}\widetilde{W}^{(N)},
\]
with all $\widetilde{K}_{-j}$ are decoupled for $0\leq j\leq N$
and ord$(\widetilde{S})=-N-1$. If we omit the term $\widetilde{S}$,
we can find an approximated solution of the form 
\[
\hat{v}_{A}^{(N)}=\sum_{j=0}^{N+1}\hat{v}_{-j,A}^{(N)}
\]
satisfying 
\[
D_{3}\hat{v}_{A}^{(N)}=\{\tau\widetilde{K_{+}^{A}}+\widetilde{K}_{0}(1,1)+\cdots+\widetilde{K}_{-N}(1,1)\}\hat{v}_{A}^{(N)}
\]
and each $\hat{v}_{-j,A}^{(N)}$ has to satisfy 
\[
\begin{cases}
D_{3}\hat{v}_{0,A}^{(N)}=\tau\widetilde{K_{+}^{A}}\hat{v}_{0,A}^{(N)}, & \hat{v}_{0,A}^{(N)}|_{x_{3}=0}=\chi_{t}(x')b,\\
D_{3}\hat{v}_{-1,A}^{(N)}=\tau\widetilde{K_{+}^{A}}\hat{v}_{-1,A}^{(N)}+\widetilde{K}_{0}(1,1)\hat{v}_{0,A}^{(N)}, & \hat{v}_{-1,A}^{(N)}|_{x_{3}=0}=0,\\
\vdots\\
D_{3}\hat{v}_{-N-1,A}^{(N)}=\tau\widetilde{K_{+}^{A}}\hat{v}_{-N-1,A}^{(N)}+\sum_{j=0}^{N}\widetilde{K}_{-j}(1,1)\hat{v}_{-j,A}^{(N)}, & \hat{v}_{-N-1,A}^{(N)}|_{x_{3}=0}=0,
\end{cases}
\]
where $\chi_{t}(x')\in C_{0}^{\infty}(\mathbb{R}^{2})$ and $b\in\mathbb{C}^{3}$.
Thus, by solving this ODE system we can get the following estimates:
\begin{equation}
\|x_{3}^{\beta}\partial_{x'}^{\alpha}(\hat{v}_{-j,A}^{(N)})\|_{L^{2}(\mathbb{R}_{+}^{3})}\leq c\tau^{-\beta-j-1/2}\label{eq:1.9}
\end{equation}
for $0\leq j\leq N+1$. Moreover, if we set $\hat{V}_{A}^{(N)}=\left[\begin{array}{c}
\hat{v}_{A}^{(N)}\\
0
\end{array}\right]$, then it satisfies 
\[
\begin{cases}
\hat{V}_{A}^{(N)}-\{\tau\widetilde{K^{A}}+\widetilde{K}_{0}+\cdots+\widetilde{K}_{-N}\}\hat{V}_{A}^{(N)}=\tilde{R},\\
\hat{V}_{A}^{(N)}|_{x_{3}=0}=\left[\begin{array}{c}
\chi_{t}(x')b\\
0
\end{array}\right],
\end{cases}
\]
where 
\[
\|\tilde{R}\|_{L^{2}(\mathbb{R}_{+}^{3})}\leq c\tau^{-N-3/2}.
\]
\textbf{Step 7.} Finally, if we define the function $\tilde{v}_{A}=\left[\begin{array}{c}
\tilde{v}_{1}\\
\tilde{v}_{2}\\
\tilde{v}_{3}
\end{array}\right]$, with $\tilde{v}_{j}$ being the $j$th component of the vector $\tilde{Q}(I+x_{3}A^{(0)}+\tau^{-1}B^{(0)})(I+x_{3}^{2}A^{(1)}+\tau^{-1}x_{3}B^{(1)}+\tau^{-2}C^{(1)})\cdots(I+x_{3}^{N+1}A^{(N)}+\tau^{-1}x_{3}^{N}B^{(N)}+\tau^{-2}x_{3}^{N-1}C^{(N)})\hat{V}_{A}^{(N)}$
and set $w_{A}=\exp(i\tau x'\cdot\xi')\tilde{v}_{A}$, we will get
that 
\begin{eqnarray*}
w_{A} & = & Q\exp(i\tau x'\cdot\xi')\exp(i\tau x_{3}\widetilde{K_{+}^{A}}(x'))\chi_{t}(x')b+\exp(i\tau x'\cdot\xi')\tilde{\Gamma}(x,\tau)\\
 & = & Q\exp(i\tau x'\cdot\xi')\exp(-i\tau x_{3}(-\widetilde{K_{+}^{A}}(x')))\chi_{t}(x')b+\Gamma(x,\tau)
\end{eqnarray*}
and 
\[
w_{A}|_{x_{3}=0}=\exp(i\tau x'\cdot\xi')(\chi_{t}(x')Qb+\beta_{0}(x',\tau),
\]
where $\beta_{0}(x',\tau)=\tilde{\Gamma}(x',0,\tau)$ is supported
in supp$(\chi_{t})$. Note that the function $\tilde{\gamma}$ comes
from the combination of $\hat{v}_{-j,A}^{(N)}$'s, for $j=1,2,\cdots,N+1$.
Now, we derive higher derivative estimates for the oscillating-decaying
solutions, back to see all the $\hat{v}_{-j,A}^{(N)}$'s separately.
In fact, only need to see $\hat{v}_{-1,A}^{(N)}$. From the estimate
(\ref{eq:1.9}), we know that the estimate is independent of the derivative
of $x'$ variables, all we need to concern is the $\partial_{3}$
derivative. From the equation 
\begin{equation}
D_{3}\hat{v}_{-1,A}^{(N)}=\tau\widetilde{K_{+}^{A}}\hat{v}_{-1,A}^{(N)}+\widetilde{K}_{0}(1,1)\hat{v}_{0,A}^{(N)}\label{eq:1.10}
\end{equation}
and the standard regularity theory of ODEs(ordinary differential equations),
we know that $\hat{v}_{-1,A}^{(N)}\in C^{\infty}$ if all the coefficients
are smooth. Moreover, note that $\widetilde{K}_{+}$ independent of
$x_{3}$, then we can differentiate (\ref{eq:1.10}) directly, to
get 
\begin{eqnarray*}
D_{3}^{2}\hat{v}_{-1,A}^{(N)} & = & D_{3}[\tau\widetilde{K_{+}^{A}}\hat{v}_{-1,A}^{(N)}+\widetilde{K}_{0}(1,1)\hat{v}_{0,A}^{(N)}]\\
 & = & \tau\widetilde{K_{+}^{A}}(D_{3}\hat{v}_{-1,A}^{(N)})+(D_{3}\widetilde{K}_{0}(1,1))\hat{v}_{0,A}^{(N)}+\widetilde{K}_{0}(1,1)D_{3}\hat{v}_{0,A}^{(N)}\\
 & = & \tau^{2}(\widetilde{K_{+}^{A}})^{2}\hat{v}_{-1,A}^{(N)}+\tau\widetilde{K_{+}^{A}}\widetilde{K}_{0}(1,1)\hat{v}_{0,A}^{(N)}+(D_{3}\widetilde{K}_{0}(1,1))\hat{v}_{0,A}^{(N)}\\
 &  & +\tau\widetilde{K}_{0}(1,1)\widetilde{K_{+}^{A}}\hat{v}_{0,A}^{(N)}.
\end{eqnarray*}
Thus, we can obtain that 
\[
\|x_{3}^{\beta}\partial_{x'}^{\alpha}\partial_{3}^{\eta}(\hat{v}_{-1,A}^{(N)})\|_{L^{2}(\mathbb{R}_{+}^{3})}\leq c\tau^{-\beta+\eta-3/2},
\]
for all $\eta\leq2$. Inductively, we have 
\[
\|x_{3}^{\beta}\partial_{x'}^{\alpha}\partial_{3}^{\eta}(\hat{v}_{-1,A}^{(N)})\|_{L^{2}(\mathbb{R}_{+}^{3})}\leq c\tau^{-\beta+\eta-3/2},
\]
for all $\eta\in\mathbb{N}$. Similarly, for other $\hat{v}_{-j,A}^{(N)}$
with $2\leq j\leq N+1$, we can get similar estimate in the following:
\[
\|x_{3}^{\beta}\partial_{x'}^{\alpha}\partial_{3}^{\eta}(\hat{v}_{-j,A}^{(N)})\|_{L^{2}(\mathbb{R}_{+}^{3})}\leq c\tau^{\eta-\beta-j-1/2}
\]
$\forall\eta\in\mathbb{N}\cup\{0\}$. Therefore, $\Gamma$ satisfies
\[
\|\partial_{x}^{\alpha}\Gamma\|_{L^{2}(\Omega_{s})}\leq c\tau^{|\alpha|-3/2}e^{-\tau(s-t)\lambda}
\]
on $\Omega_{s}:=\{x_{3}>s\}\cap\Omega$ for $s\geq0$ and $\forall|\alpha|\in\mathbb{N}\cup\{0\}$.
Note that since each $\hat{v}_{-j,A}^{(N)}$'s are smooth, we can
get the smoothness of $\tilde{R}$ and 
\[
\|\partial_{x}^{\alpha}\tilde{R}\|_{L^{2}(\mathbb{R}_{+}^{3})}\leq c\tau^{|\alpha|-N-3/2}
\]
for all $|\alpha|\in\mathbb{N}\cup\{0\}$. Furthermore, we have that
\[
\|\partial_{x}^{\alpha}(Q_{A}\widetilde{v_{A}})\|_{L^{2}(\Omega_{0})}\leq c\tau^{|\alpha|-N-1/2}.
\]
\textbf{Step 8.} Now let $u=w+r=e^{i\tau x'\cdot\xi'}\tilde{v}+r$
and $r$ be the solution to the boundary value problem 
\[
\begin{cases}
L_{A}r=-e^{i\tau x'\cdot\xi'}\widetilde{Q_{A}}\widetilde{v_{A}} & \mbox{ in }\Omega_{0}\\
r=0 & \mbox{ on }\partial\Omega_{0}
\end{cases}.
\]
However, note that $\Omega_{0}=\{x_{3}>0\}\cap\Omega$ is not a smooth
domain since $\partial\Omega_{0}=(\{x_{3}=0\}\cap\Omega)\cup(\{x_{3}>0\}\cap\partial\Omega)$.
Note that the oscillating-decaying solution exists in the half space,
from the construction, we know that the solution is independent of
the domain $\Omega$. Let $\widetilde{\Omega}\subset\mathbb{R}_{+}^{3}$
be a open bounded smooth domain containing $\Omega$ with $\{x_{3}=0\}\cap\Omega\subset\partial\widetilde{\Omega}$,
from the construction, it is easy to see the form of oscillating-decaying
solution does not depend on the domain $\Omega$, then we can extend
$r$ to be defined on $\widetilde{\Omega}$ and call it $\tilde{r}(x)$.
Here we can also extend $\widetilde{v_{A}}$ to be defined on $\widetilde{\Omega}$,
still denote $\widetilde{v_{A}}$ and all the decaying estimates will
hold since our estimates were considered in $\mathbb{R}_{+}^{3}$,
then we have 
\[
\begin{cases}
L_{A}\tilde{r}=-e^{i\tau x'\cdot\xi'}\widetilde{Q_{A}}\widetilde{v_{A}} & \mbox{ in }\widetilde{\Omega},\\
\tilde{r}=0 & \mbox{ on }\partial\widetilde{\Omega}.
\end{cases}
\]
Note that all the coefficients are smooth, we apply a well-known elliptic
regularity theorem (Theorem2.3, \cite{G1993book}), then we will get
$\tilde{r}\in C^{k}(\Omega)$ $\forall k$ (recall that $\partial\Omega\in C^{\infty}$)
and 
\[
\|\tilde{r}\|_{H^{k+1}(\Omega;\mathbb{R}^{3})}\leq c\|\widetilde{Q_{A}}\widetilde{v_{A}}\|_{H^{k}(\Omega;\mathbb{R}^{3})}.
\]
Hence $\|\partial_{x}^{\alpha}r\|_{L^{2}(\Omega_{0})}\leq\|\partial_{x}^{\alpha}\tilde{r}\|_{L^{2}(\widetilde{\Omega})}\leq c\tau^{|\alpha|-N+1/2}$
for all $|\alpha|\leq k$, $\forall k\in\mathbb{N}$. Similarly, we
can construct the oscillating decaying solution for $L_{B}B=0$. Then
we represent $A$ and $B$ to be two oscillating-decaying solution
in the following form: 
\[
\begin{cases}
A= & w_{\chi_{t},b,t,N,\omega}^{A}+r_{\chi_{t},b,t,N,\omega}^{A}\mbox{ in }\Omega_{t}(\omega),\\
A= & e^{i\tau x\cdot\xi}\{\chi_{t}(x')Q_{t}(x')b+\beta_{\chi_{t},t,b,N,\omega}^{A}\}\mbox{ on }\Sigma_{t}(\omega),\\
B= & w_{\chi_{t},b,t,N,\omega}^{B}+r_{\chi_{t},b,t,N,\omega}^{B}\mbox{ in }\Omega_{t}(\omega),\\
B= & e^{i\tau x\cdot\xi}\{\chi_{t}(x')Q_{t}(x')b+\beta_{\chi_{t},t,b,N,\omega}^{B}\}\mbox{ on }\Sigma_{t}(\omega),
\end{cases}
\]
where 
\[
\begin{cases}
w_{\chi_{t},b,t,N,\omega}^{A}=\chi_{t}(x')Q_{t}e^{i\tau x\cdot\xi}e^{-\tau(x\cdot\omega-t)A_{t}^{A}(x')}b+\gamma_{\chi_{t},b,t,N,\omega}^{A}(x,\tau),\\
w_{\chi_{t},b,t,N,\omega}^{B}=\chi_{t}(x')Q_{t}e^{i\tau x\cdot\xi}e^{-\tau(x\cdot\omega-t)A_{t}^{B}(x')}b+\gamma_{\chi_{t},b,t,N,\omega}^{B}(x,\tau),
\end{cases}
\]
$\gamma_{\chi_{t},b,t,N,\omega}^{A}$ and $\gamma_{\chi_{t},b,t,N,\omega}^{B}$
satisfy (\ref{eq:1.6}) and (\ref{eq:1.7}).

\subsection{Well-posedness and $L^{p}$ estimate for the anisotropic Maxwell
system}

In the following, we would list the eigenvalue property and well-posedness
results of the following problem: let $\Oa\subset\mb R^{3}$ and $K\Subset\Oa$,

\begin{equation}
\begin{cases}
\na\times E=ik\mu H & \quad\mbox{in }\Oa\setminus K\\
\na\times H=-ik\epsilon E+J & \quad\mbox{in }\Oa\setminus K\\
\nu\times E=f & \quad\mbox{on }\pl\Oa\\
\nu\times H=g & \quad\mbox{on }\pl K,
\end{cases}\label{eq:Pb}
\end{equation}
where $\mu,\en$ are symmetric and positive definite matrix-valued
functions. More precisely, we assume there exist constants $\mu_{0},\mu_{1},\la_{0},\La_{0}>0$
such that 
\begin{equation}
\begin{cases}
\mu_{0}I\leq\mu(x)\leq\mu_{1}I,\\
\la_{0}I\leq\epsilon(x)\leq\La_{0}I.
\end{cases}\label{eq:assump}
\end{equation}
These well-posedness for the isotropic Maxwell systems can be found
in Theorem 4.18 and 4.19 of \cite{M2003book}. However, we have the
same result under our assumption (\ref{eq:assump}) following the
arguments in \cite{M2003book}. Let 
\[
X=\left\{ u\in H(curl;\Oa\setminus K)|\nu\times u=0\mbox{ on }\pl\Oa\mbox{ and }u_{T}\in L^{2}\left(\pl K\right)^{3}\mbox{ on }\pl K\right\} .
\]

\begin{defn}
We say $(E,H)$ or $E$ is a weak solution of (\ref{eq:Pb}) if $E\in X$
and satisfies
\begin{equation}
\left\langle \mu^{-1}\na\times E,\na\times\phi\right\rangle _{\Oa\setminus K}-k^{2}\left\langle \gamma E,\phi\right\rangle _{\Oa\setminus K}=\left\langle ikJ,\phi\right\rangle _{\Oa\setminus K}-\left\langle \mu^{-1}g,\phi_{T}\right\rangle _{\pl K},\quad\forall\phi\in X,\label{eq:WF}
\end{equation}
and $\bs{\nu}\times E=f$ on $\pl\Oa$, where $\phi_{T}=\left(\nu\times\phi\right)\times\nu$
and $\left\langle \cdot,\cdot\right\rangle $ denotes the standard
Hermitian inner product of $L^{2}$ space. Moreover, if (\ref{eq:WF})
fails to have a unique solution, then $k$ is called an eigenvalue
or a resonance of (\ref{eq:Pb}). \end{defn}
\begin{lem}
\label{eigen prop}There is an infinite discrete set $\Sigma$ of
eigenvalue $k_{j}>0$, $j=1,2,\ldots$ and corresponding eigenfunctions
$E_{j}\in H_{0}(curl;\Oa)$, $E_{j}\neq0$, such that (\ref{eq:WF})
holds with $J=0$ and $f=g=0$ is satisfied.
\end{lem}
From the above lemma, we have the following theorem. 
\begin{thm}
\label{well pose} For $k\notin\Sigma$, there exists a unique weak
solution $(E,H)\in H(curl;\Omega\backslash\ol K)\times H(curl;\Omega\backslash\ol K)$
of (\ref{eq:Pb}) given any $f\in H^{-1/2}(Div;\pl\Omega)$, $g\in H^{-1/2}(Div;\partial K)$
and $J\in H^{-1}(\Omega\backslash\ol K)$. The solution satisfies
\[
\|E\|_{L^{2}(\Omega\backslash\ol K)}+\|H\|_{L^{2}\left(\Oa\setminus K\right)}\leq C(\|f\|_{H^{-1/2}(Div;\partial\Omega)}+\|g\|_{H^{-1/2}(Div;\partial K)}+\|J\|_{H^{-1}(\Omega\backslash\ol K)})
\]
for some constant $C>0$, where
\[
H^{-1/2}\left(Div;\Ga\right):=\left\{ f\in H^{-1/2}\left(\Ga\right)^{3}\left|\begin{array}{c}
\nu\cdot f=0,\mbox{ }\na_{\pl\Oa}\cdot f\in H^{-1/2}\left(\Ga\right)\end{array}\right.\right\} ,
\]
$\Ga=\pl\Oa$ or $\pl K$.
\end{thm}
In the following, we state the $L^{p}$ theory for the anisotropic
Maxwell's system. For this purpose, we define a bilinear form 
\[
B_{A}(E,F):=\int_{\Omega}(A(x)\nabla\times E(x))\cdot(\nabla\times\overline{F}(x))dx+M\int_{\Omega}E(x)\cdot\overline{F}(x)dx
\]
for all $E\in H_{0}^{1,q}(curl,\Omega)$ and $F\in H_{0}^{1,q'}(curl,\Omega)$
with $\dfrac{1}{q}+\dfrac{1}{q'}=1$. We only state $L^{p}$ estimate
in the following theorem, but we do not prove the theorem. For more
details, we refer readers to read \cite{KS2014}.
\begin{thm}
\cite{KS2014} Let $\Omega$ be a a smooth domain. Suppose that $A=A(x)$
is a real symmetric matrix with smooth entries and satisfies the uniform
elliptic condition 
\[
\lambda|\xi|^{2}\leq A(x)\xi\cdot\xi\leq\Lambda|\xi|^{2},\mbox{ for all }\xi\in\mathbb{R}^{3},
\]
for some constants $0<\lambda\leq\Lambda<\infty$. Assume $q$ is
some number satisfying $2\leq q<\infty$. Under the condition
\[
\inf_{\|F\|_{1,q'}=1}\sup_{\|E\|_{1,q}=1}|B_{A}E,F)|\geq\dfrac{1}{K}>0
\]
the Maxwell's systems of the equations 
\[
\nabla\times(A\nabla\times E)+E=\nabla\times f+g
\]
is uniquely solvable in $H_{0}^{1,q'}(curl,\Omega)$ for each $g\in L^{q'}(\Omega)$
and $f\in L^{q'}(\Omega)$ and the weak solution satisfies 
\[
\|E\|_{L^{q'}(\Omega)}+\|\nabla\times E\|_{L^{q'}(\Omega)}\leq K\{\|f\|_{L^{q'}(\Omega)}+\|g\|_{L^{q'}(\Omega)}\},
\]
where $K$ is a positive constant depending on $p$.
\end{thm}
We end up this appendix with the following lemma on the embedding
related to the Sobolev-Besov spaces, for more details, see \cite{mitrea1996vector}.
\begin{lem}
Let $u\in L^{p}(D)$ such that $\nabla\cdot u\in L^{p}(D)$ and $\nabla\times u\in L^{p}(D)$.
If $\nu\times u\in L^{p}(\partial D)$, then also $\nu\cdot u\in L^{p}(\partial D)$
for $p\in(1,\infty)$. If in addition $1<p\leq2$, then $u\in B_{\frac{1}{p}}^{p,2}(D)$
and we have the estimate 
\[
\|u\|_{B_{\frac{1}{p}}^{p,2}(D)}\leq C\{\|u\|_{L^{p}(D)}+\|\mbox{curl}u\|_{L^{p}(D)}+\|\nabla\cdot u\|_{L^{p}(D)}+\|\nu\times u\|_{L^{p}(\partial D)}\}
\]
where the Sobolev-Besov space $B_{\alpha}^{p,q}(D):=[L^{p}(D),W^{1,p}(D)]_{\alpha,q}$
is obtained by real interpolation for $1<p,q<\infty$ and $0<\alpha<1$.
\end{lem}
\bibliographystyle{plain}
\bibliography{ref}

\end{document}